\documentclass[11pt]{amsart}
\usepackage[cal=cm]{mathalfa}
\usepackage{pgfplots, graphicx, subcaption, amsmath, array, pdfsync, tikz, tikz-cd, color, url, float, fullpage, wrapfig,comment, scrextend, wrapfig, amsthm, colonequals,amssymb, enumitem}
\usepackage[normalem]{ulem}

\DeclareMathOperator{\Cliff}{\mathrm{Cliff}}

\newcommand{\floor}[1]{\left\lfloor #1 \right\rfloor}
\newcommand{\ol}[1]{\overline{#1}}

\newcommand{\g}{\mathfrak{g}}
\newcommand{\ot}{\otimes}
\DeclareMathOperator{\Pol}{\mathrm{Pol}}
\DeclareMathOperator{\im}{\mathrm{im}}
\newcommand{\embed}{\hookrightarrow}

\usetikzlibrary{decorations.pathreplacing,calligraphy}

\usepackage[pagebackref,pdfencoding=auto]{hyperref}
\hypersetup{
  colorlinks= true,
  linkcolor=[RGB]{0,89,42},
  urlcolor=[RGB]{0,89,42},
  citecolor=[RGB]{0,89,42}
}
\usepackage{cleveref}

\usepackage[font=footnotesize]{caption}

\usetikzlibrary{calc}

\title{Irrationality and Noether--Lefschetz loci in the moduli space of K3 surfaces}
\author{Dori Bejleri, Nathan Chen, Benjamin Church, Kristin DeVleming}

\definecolor{mycolor}{RGB}{146, 214, 203}
\definecolor{myothercolor}{RGB}{179, 215, 232}

\definecolor{mycolor}{RGB}{146, 214, 203}
\definecolor{myothercolor}{RGB}{179, 215, 232}

\newtheorem{theorem}{Theorem}
\newtheorem{corollary}[theorem]{Corollary}
\newtheorem{lemma}[theorem]{Lemma}
\newtheorem{proposition}[theorem]{Proposition}
\newtheorem{prop}[theorem]{Proposition}

\numberwithin{theorem}{section}

\theoremstyle{definition}
\newtheorem{remark}[theorem]{Remark}

\newtheorem{example}[theorem]{Example}
\newtheorem{definition}[theorem]{Definition}

\renewcommand{\P}{\mathbb{P}}
\def\C{{\mathbb C}}
\def\F{{\mathbb F}}

\def\Zb{{\mathbb Z}}

\def\Lc{{\mathcal L}}
\def\Fc{{\mathcal F}}
\def\Gc{{\mathcal G}}
\def\Oc{{\mathcal O}}

\def\Sc{{\mathcal S}}
\def\Xc{{\mathcal X}}

\def\Qb{{\mathbb Q}}

\def\Ac{{\mathcal A}}
\def\Hc{{\mathcal H}}

\def\NL{{\mathrm{BN}}}
\newcommand{\gonlocus}[2]{\mathcal{G}_{#1 ; #2}}
\newcommand{\pgonlocus}[2]{\mathcal{P}_{#1 ; #2}}
\DeclareMathOperator{\Cliffdim}{\mathrm{Cliff.dim}}

\def\irr{{\mathrm{irr}}}
\def\NL{{\mathrm{NL}}}

\def\Aut{{\mathrm{Aut}}}
\def\Bir{{\mathrm{Bir}}}

\DeclareMathOperator{\Sym}{Sym}

\DeclareMathOperator{\GL}{GL}

\DeclareMathOperator{\Pic}{\mathrm{Pic}}
\DeclareMathOperator{\gon}{\mathrm{gon}}

\DeclareMathOperator{\sat}{\mathrm{sat}}

\definecolor{greener}{RGB}{72,171,131}

\pgfplotsset{compat=1.15}

\begin{document}

\begin{abstract}
    The goal of this paper is twofold. First, we construct examples of families of surfaces where the locus of points in the base with constant degree of irrationality is not constructible. Second, we explore the notion of gonality loci in the moduli space of K3 surfaces, which are characterized by the existence of a smooth curve with fixed gonality. Among various results, we prove that polarized gonality loci are always unions of Noether--Lefschetz divisors.
\end{abstract}

\maketitle

\section{Introduction}

A fundamental result in rationality is Castelnuovo's criterion (due to Zariski in positive characteristic), which states that a smooth projective surface $S$ over an algebraically closed field is rational if and only if $h^{1}(\Oc_S) = 0$ and $h^{0}(\omega_{S}^{\otimes 2}) = 0$. As a consequence, if $\Sc \rightarrow T$ is a family of integral surfaces over $\C$ then the locus of closed points $t \in T$ such that $\Sc_{t}$ is rational is a constructible set. In higher dimensions, Hassett--Pirutka--Tschinkel \cite{HPT18} showed that a cohomological criterion for rationality does not exist, by constructing families $\Xc \rightarrow B$ of quadric surface bundles such that (1) for very general $b \in B$, the fiber $\Xc_{b}$ is not stably rational, but (2) the set of points $b \in B$ where $\Xc_{b}$ is rational is dense in $B$ for the Euclidean topology.

Motivated by the results above, the first part of this paper is devoted to studying irrationality loci, which are subsets of the base of a family where the fibers have constant degree of irrationality. Recall that the degree of irrationality $\irr(X)$ of an $n$-dimensional variety $X$ is by definition the least degree of a rational covering $X \dashrightarrow \P^{n}$. We construct examples showing that irrationality loci are not constructible even in smooth families of surfaces, and state this in the same spirit as \cite{HPT18}.

\begin{theorem}\label{thm:main-result}
There exist smooth families of complex projective surfaces $\Sc \rightarrow T$ over connected varieties $T$ such that
\begin{enumerate}
    \item for very general $t \in T$, the fiber $\Sc_{t}$ satisfies $\irr(\Sc_{t}) = d$ for some integer $d \geq 3$, and
    \item the set of points $t \in T$ such that $\Sc_{t}$ satisfies $\irr(\Sc_{t}) < d$ is dense in $T$ for the Euclidean topology.
\end{enumerate}
\end{theorem}

\noindent The condition that $d \geq 3$ is necessary, since the same statement for $d = 2$ is false by Castelnuovo's criterion. In fact, \cite[Corollary 3.1]{CE25} implies that in a smooth family, having $\irr(S_{t}) = 2$ for $t \in T$ very general implies the same statement for all $t \in T$.

The theorem above relies on an analysis of the locus of K3 surfaces $X$ with $\irr(X) = 2$, for which we establish two lattice-theoretic criteria (see \S\ref{section:hyperellipticK3s}). Along the way, we study certain ``gonality loci" in the moduli space $\Fc_{g}$ of quasi-polarized K3 surfaces $(X, L)$ of genus $g \geq 2$, which are defined as follows.
\[ \gonlocus{e}{g} \colonequals \{ (X, L) \in \Fc_{g} \colon \exists \text{ smooth curve } C \subset X \text{ with } \gon(C) = e \}. \]
Note that there is a priori no connection between $L$ and $C$, so we refer to these as \textit{unpolarized} gonality loci. Our analysis of these loci for $e = 2$ gives a full description of $\gonlocus{2}{g}$ and its interaction with the locus $\Hc_{g} \subset \Fc_{g}$ of K3 surfaces that have $\irr = 2$. For instance, from \Cref{thm:irr2hypcurve} we will see that $\Hc_{g} \subset \gonlocus{2}{g}$ and coincides with the locus of pairs $(X, L) \in \gonlocus{2}{g}$ for which the curve $C \subset X$ with $\gon(C) = 2$ has genus at least two. 

\begin{theorem}\label{thm:Hg-in-G2g}
Fix $g \geq 3$. If $g$ is even, then $\Hc_{g} \subsetneq \gonlocus{2}{g}$ consists of a union of countably infinitely many divisors. If $g$ is odd, then $\Hc_{g} \subsetneq \gonlocus{2}{g}$ consists of a union of countably infinitely many divisors and countably infinitely many codimension 2 loci. 
\end{theorem}

\begin{corollary}\label{cor:G2g}
For $g \geq 3$, the locus $\gonlocus{2}{g}$ is a countably infinite union of divisors in $\Fc_{g}$.
\end{corollary}

\noindent When $g = 2$, it will follow from \Cref{lem:K3DoubleCoverP2} and \Cref{prop:base-locus-of-big-and-nef} that $\Fc_{2} = \gonlocus{2}{2}$. For $g \geq 3$ and $e = 2$, see \Cref{thm:Irr2K3s}, \Cref{thm:Irr2K3s-version2}, and \Cref{thm:codim2-hyperelliptic} for a more precise version of \Cref{thm:Hg-in-G2g}. One of the original goals of this project was to understand the relationship between the degree of irrationality of a K3 surface and the existence of smooth higher gonality curves, but this is much less clear (see \S\ref{section:irr3-and-trigonal} for a relevant example).

The last section is devoted to understanding \textit{polarized} gonality loci $\pgonlocus{e}{g}^{m}$, where the curve $C \subset X$ of gonality $e$ is required to be linearly equivalent to a fixed multiple $mL$ of the quasi-polarization. In the interesting cases, we show that they are always divisorial.

\begin{theorem}\label{intro-thm:pol-gonloci-divisors}
For fixed integers $e,g,m > 0$, the locus $\ol{\pgonlocus{e}{g}^{m}}$ is either empty, a finite union of Noether--Lefschetz divisors in $\Fc_g$, or all of $\Fc_g$.
\end{theorem}

Finally, in the Appendix we give a checkable criterion for the existence of a codimension $\geq 2$ component of $\gonlocus{e}{g}$ in terms of a certain $3 \times 3$ sublattice of the K3 lattice (Theorem \ref{thm:codim-2-comp}), and use this to show that unpolarized gonality loci need not be divisorial when $e \ge 3$.  This generalizes the behavior of $\Hc_g$ in the odd genus case of Theorem \ref{thm:Hg-in-G2g}.

\begin{theorem}\label{intro-thm-codim2}
    For $e \ge 3$, there exist components of $\gonlocus{e}{g}$ that have codimension 2 and are not purely divisorial.  
\end{theorem}

\subsection*{Outline} In \S\ref{section:prelims}, we collect some well-known results about K3 surfaces as well as the Clifford index and gonality of curves in a given linear series. In \S\ref{section:Irr2HyperellipticCurves}, we describe some consequences of a surface $X$ having $\irr(X) = 2$ and then specialize to the case of K3 surfaces. In \S\ref{section:hyperellipticK3s} we give two equivalent formulations for a K3 surface to have $\irr(X) = 2$, and use these to prove \Cref{thm:Hg-in-G2g}. In \S\ref{section:IrratLociNotConstructible} we prove \Cref{cor:G2g} and \Cref{thm:main-result}. In \S\ref{section:irr3-and-trigonal}, we present an example of a K3 surface with degree of irrationality 3 that does not contain a smooth trigonal curve. In the last section, we establish some properties about polarized gonality loci and prove \Cref{intro-thm:pol-gonloci-divisors}.  Theorem \ref{intro-thm-codim2} is proved in the appendix.

\subsection*{Acknowledgments} We would like to thank R. Haburcak, B. Hassett, R. Lazarsfeld, O. Martin, J.C. Ottem, and G. P. Pirola for stimulating discussions. After an early draft of the paper was written up, NC learned that J. C. Ottem had outlined a similar approach to \Cref{thm:main-result}. DB was supported by NSF grant DMS-2401483 and NSF CAREER grant DMS-2542956. KD was supported by NSF grant DMS-2550445 and NSF CAREER grant DMS-2542575. We would also like to thank the Brin Center and the University of Maryland for their hospitality during a week-long research visit.

\subsection*{AI Disclosure}

In the preparation of this paper, AI models were used to find examples of matrices that satisfy the assumptions of \Cref{thm:codim-2-comp}. They were also used to audit the paper after it was completed, and found a missing case in the original proof of Theorem~\ref{thm:polLociDivisors}. The authors have verified the details and take full responsibility for the correctness.

\section{Preliminaries}\label{section:prelims}

Throughout this paper, we will work over $\C$.

\begin{definition}
    A \textit{K3 surface} is a smooth projective surface $X$ defined over $\C$ such that $\omega_X \cong \Oc_X$ and $h^1(X, \Oc_X) = 0$. A \textit{(quasi-)polarized K3 surface} is a pair $(X,L)$ such that $X$ is a K3 surface and $L$ is a primitive ample (resp. big and nef) line bundle on $X$. The \textit{degree} of a quasi-polarized K3 surface $(X,L)$ is $d \colonequals L^2 > 0$ and the \textit{genus} $g$ of $X$ satisfies $d = 2g -2$.
\end{definition}

If the complete linear series $|L|$ is basepoint-free, then $g$ is also the genus of a smooth curve in $|L|$. The following fact about K3 surfaces will be used frequently.
 
\begin{proposition}[{\cite{Huybrechts16, Mayer, Saint-Donat74}}]\label{prop:base-locus-of-big-and-nef}
    If $L$ is a line bundle on a K3 surface $X$, then $|L|$ has no base points outside its fixed locus. If $L$ is in addition big and nef, then:
        \begin{enumerate}
            \item if the mobile part $M$ of the linear series $|L|$ is big, then $L$ is basepoint-free and $L = M$;
            \item if $L$ is not basepoint-free, then $L \cong \Oc(mE + R)$ where $E \subset X$ is a smooth elliptic curve, $R \subset X$ is a smooth rational curve, $E \cdot R = 1$, and $m \ge 2$. Moreover, $eL$ is basepoint-free for $e \geq 2$.
        \end{enumerate}
\end{proposition}

\noindent This appears in several places in slightly different forms. \cite[Chapter 2, Corollary 3.15 \& Remark 3.4]{Huybrechts16} is almost the same statement. In part (2), the fact that $E \cdot R = 1$ follows from \cite[Proposition 5]{Mayer} since we know that $R$ consists of a single rational curve, and the statement that $eL$ is basepoint-free for $e \geq 2$ follows from \cite[Chapter 2, Remark 3.4]{Huybrechts16}

We will also make use of \textit{reflections} as defined in the following remark, which allow us to assume certain line bundles are nef. 

\begin{remark}\label{rem:reflections-nef}
    Let $L$ be a line bundle on a K3 surface $X$ with $L^2 \ge 0$.  Fix an ample line bundle $H$ on $X$.  Up to replacing $L$ with its dual, by Riemann--Roch, we may assume $L$ is effective.  If $L$ is not nef, there exists a rational curve $R \subset X$ such that $L \cdot R < 0$.  The \textit{reflection} of $L$ across the $(-2)$-curve $R$ is the line bundle $L' = L +(L\cdot R) R$.  Observe that $(L')^2 = L^2 + 2(L \cdot R)^2 - 2(L \cdot R)^2 = L^2$ but $L' \cdot H < L \cdot H$.  If $L'$ is not nef, we repeat, and as the intersection with the fixed ample class $H$ decreases after each reflection, this process terminates after finitely many steps. Therefore, if $L$ is a line bundle on a K3 surface with $L^2 \ge 0$, then there exists a nef line bundle $L'$ such that $L'^2 = L^2$.
\end{remark}

\begin{definition}
The Noether--Lefschetz divisor $\NL_{\Lambda;g} \subset \Fc_g$ is the locus of polarized K3 surfaces $(X,L)$ such that $L \in \Lambda \subset \Pic{X}$ where $\Lambda \subset \Lambda_{K3}$ is a primitively embedded sublattice of the K3 lattice of rank $2$ and signature $(1,1)$.    
\end{definition}

\subsection{Lifting divisors on curves in a linear series}

The influential work of Saint-Donat, Green--Lazarsfeld, Donagi--Morrison, Knutsen, and many others shows that we can completely compute many invariants of a smooth curve on a K3 surface lattice-theoretically. We begin by recalling some definitions.

\begin{definition}
    Let $C$ be a smooth curve of genus $g \ge 2$ and let $L$ be a line bundle on $C$. Recall that the \textit{Clifford index} of $L$ on $C$ is 
    \[ \Cliff(L) = \deg L - 2(h^0(L) - 1)\] and the \textit{Clifford index} of $C$ is 
    \[ \Cliff(C) = \min \{ \Cliff(L) \mid h^0(L) \ge 2, h^1(L) \ge 2 \}.\]
    We say that $L$ \textit{contributes to the Clifford index} if $h^0(L) \ge 2$ and $h^1(L) \ge 2$ and we say that $L$ \textit{computes the Clifford index} of $C$ if $\Cliff(L) = \Cliff(C)$. The \textit{Clifford dimension} (denoted $\Cliffdim(C)$) is the minimal $r$ so that the Clifford index is computed by a $\g^r_d$ on $C$.
\end{definition}

We say that the Clifford index is \textit{generic} if it is equal to its maximum (and generic) value $\floor{\frac{g-1}{2}}$. We also recall that a curve $C$ is \textit{exceptional} if $\Cliffdim(C) > 1$. Concretely, $C$ is exceptional if and only if $\Cliff(C) = \gon(C) - 3$, whereas $\Cliff(C) = \gon(C) - 2$ when $\Cliffdim(C) = 1$. 

The following theorem of Green and Lazarsfeld \cite{GL87} will play a key role.

\begin{theorem}[{\cite{GL87}}]\label{thm:GL}
    Let $X$ be a K3 surface and $C \subset X$ be a smooth curve of genus $g \ge 2$.  Then:
    \begin{enumerate}
        \item For any smooth curve $C' \in |C|$, we have $\Cliff(C') = \Cliff(C)$. 
        \item If $C$ is a smooth curve on $X$ such that $\Cliff(C) $ is non-generic, i.e. $\Cliff(C) < \floor{\frac{g-1}{2}}$, then there is a line bundle $L$ on $X$ such that $\Cliff(C')$ is computed by the restriction of $L$ to $C'$.
    \end{enumerate}
\end{theorem}

In contrast, a $\g^1_d$ realizing the gonality of the curve need not come from the K3 surface. However, a result of Donagi--Morrison roughly says that the gonality pencil on the curve $C$ is contained in the restriction of an effective divisor on the K3 surface with controlled Clifford index. 

\begin{theorem}[{\cite[Theorem 5.1]{DonagiMorrison}}]\label{thm:DonagiMorrison}
Let $X$ be a K3 surface and $C \subset X$ a smooth nonhyperelliptic curve of genus $g$ carrying a complete, basepoint-free $\g^1_d$ denoted by $|Z|$ whose corresponding Brill--Noether number satisfies $\rho(g,1,d) < 0$. Then there exist effective divisors $D, \Delta$ on $X$ such that 
\begin{enumerate}[label=(\alph*)]
\item $D + \Delta \sim C$
\item $h^0(X, \Oc_X(D)) \ge 2$ and $h^0(X, \Oc_X(\Delta)) \ge 2$
\item $\Cliff(D|_C) \le \Cliff(Z)$
\item $\Cliff(D|_{C'}) = \Cliff(D|_C)$ for all $C' \in |C|$ in the linear series
\item there exists a divisor $D_0 \in |D|$ and a reduced divisor $Z_0 \in |Z|$ such that $Z_0 \subset D_0 \cdot C$
\item $d \le D \cdot C \le g - 1$.
\item $2 \le D \cdot \Delta \le d$
\end{enumerate}
\end{theorem}

\noindent The first inequality in part (f) follows from part (e), while part (g) follows from \cite[Corollary~4.5]{DonagiMorrison} and the proof of \cite[Lemma~4.6]{DonagiMorrison} which shows that $0 \le \Cliff(D|_C) \le D \cdot \Delta - 2$.

\begin{remark}\label{rmk:DM-g1d}
    If the $\g^1_d$ in the previous theorem matches up with $\gon(C) = d$ and $C$ is not exceptional, then in fact $D \cdot \Delta = d$. To see this, note that $d = \gon(C) = \Cliff(C) + 2 \le \Cliff(D|_C) + 2 \le D\cdot \Delta$, and then combine with part (g) in the previous theorem.  If $C$ is exceptional, then $d = \gon(C) = \Cliff(C) + 3 \le \Cliff(D|_C) + 3 \le D \cdot \Delta + 1$ combined with part (g) implies $\Cliff(C)$ is at most one less than $\Cliff(D|_C)$, i.e. $\Cliff(D|_C) - 1 \le \Cliff(C) \le \Cliff(D|_C)$. 
\end{remark}

\begin{remark}
As stated, the theorem of Donagi--Morrison only applies to \textit{complete} $\g^1_d$ linear series. This is harmless for our purposes since we only apply the theorem to extremal $\g^1_d$ computing the gonality. Indeed, if $f : C \to \P^1$ is a map of minimal degree $d = \gon(C)$, then $A =  f^* \Oc_{\P^1}(1)$ necessarily satisfies $h^0(C,A) = 2$. Otherwise, if $h^0(C, A) \ge 3$, then for any $p \in C$ the line bundle $A(-p)$ has degree $d-1$ and $h^0(C, A(-p)) \ge 2$. After removing base points, the new line bundle defines a map $C \to \P^1$ of degree $\le d - 1$, a contradiction.
\end{remark}

Recall that the Brill--Noether number $\rho(g,1,d) = 2d - g - 2$ is negative if and only if $d < \floor{\frac{g+3}{2}}$, so the theorem above applies to curves with non-generic gonality for their genus. Our main use of this theorem is to show that the existence of a curve of non-generic gonality forces one of a finite list (depending on the genus and gonality of the curve) of $2 \times 2$ sublattices to appear in $\Pic{X}$. 

\begin{definition}\label{def:DMlattice}
Let $(g,e)$ be a pair of positive integers. A \textit{Donagi--Morrison matrix} of type $(g, e)$ is a symmetric $2 \times 2$ integer matrix of the form
\[ \begin{pmatrix}
a & b
\\
b & c
\end{pmatrix}\]
where $a,b,c \in \Zb$ satisfy:
\begin{enumerate}[label=(\alph*)]
    \item $a$ and $c$ are even integers
    \item $0 \le a \le c$
    \item $a + 2 b + c = 2 g - 2$
    \item $2 \le b \le e$
    \item $ac \le b^2$
\end{enumerate}
\end{definition}

For each pair $(g, e)$ of positive integers, observe that there are only finitely many Donagi--Morrison matrices of type $(g,e)$.  Examples of these matrices will arise as the intersection matrices of $D$ and $\Delta$ as in the previous theorem.

\begin{corollary}\label{cor:DMmatrices}
Let $X$ be a K3 surface and $C \subset X$ a smooth nonhyperelliptic curve of genus $g$ carrying $|Z|$ a complete, basepoint-free $\g^1_d$ whose corresponding Brill--Noether number satisfies $\rho(g,1,d) < 0$. Then $\Pic{X}$ contains 2 classes $D, \Delta$ whose span contains $C$ and whose intersection form is represented by a $(g, d)$ Donagi--Morrison matrix.
\end{corollary}

\begin{proof}
Let $D, \Delta$ be as in the conclusion of Theorem \ref{thm:DonagiMorrison} and set $\Lambda = \left< D, \Delta \right>$ . Then $C \sim D + \Delta \in \Lambda$ by Theorem \ref{thm:DonagiMorrison}(a). Moreover, the intersection matrix in the $D, \Delta$ basis for $\Lambda$ is a Donagi--Morrison matrix.  Indeed, define:
\[ \begin{pmatrix}
a & b
\\
b & c
\end{pmatrix}
=
\begin{pmatrix}
D^2 & D \cdot \Delta 
\\
D \cdot \Delta & \Delta^2
\end{pmatrix}
\]
From Theorem \ref{thm:DonagiMorrison}(f) and (g), we have $a$ and $c$ are non-negative. Furthermore, 
\[ 2 g - 2 = C^2 = D^2 + 2 D \cdot \Delta + \Delta^2 = a + 2b + c. \]
Now from Theorem \ref{thm:DonagiMorrison}(f), we have $d \le D \cdot C \le g - 1$ and hence
we obtain $d \le a + b \le g - 1$. Therefore, using the fact that $2g -2 = a+2b+c$, we obtain $a - c \le 0$ or $a \le c$. Next, from Theorem \ref{thm:DonagiMorrison}(g) we have the inequality $2 \le b = D \cdot \Delta \le d$. Finally the Hodge index theorem implies the determinantal condition $ac \le b^2$, with equality holding only in the case that $D$ and $\Delta$ are proportional.
\end{proof}

Note that the existence of a $(g,e)$ Donagi--Morrison matrix in $\Pic{X}$ forces the existence of curves with bounded Clifford index and gonality.

\begin{lemma}\label{lem:DM-matrices-force-gon(C)}
    Let $X$ be a K3 surface and suppose $M = \begin{pmatrix} a & b \\ b & c \end{pmatrix}$ is a $2 \times 2$ Donagi--Morrison matrix of type $(g,e)$ spanned by classes $D$ and $\Delta$ in $\Pic{X}$.  Then, we may assume $D$ and $\Delta$ are effective with at least two sections and $D+\Delta$ is big and nef and either 
    \begin{enumerate}
        \item $D+ \Delta$ is not basepoint-free, in which case the Picard lattice of $X$ represents $0$ and $-2$, or 
        \item the linear series $|D + \Delta|$ contains a smooth irreducible curve $C$ of Clifford index at most $b -2$.  In particular, if $C$ is not exceptional, $\gon(C) \le b$, and if $C$ is exceptional, $\gon(C) \le b+1$. 
    \end{enumerate}
\end{lemma}

\begin{proof}
    Let $M = \begin{pmatrix} a & b \\ b & c \end{pmatrix}$ satisfying Definition \ref{def:DMlattice} and denote by $D, \Delta$ the classes spanning $M$.  Because $0 \le a \le c$, we have $0 \le D^2$ so either $D$ or $-D$ is effective. Up to replacing $D, \Delta$ by $-D, - \Delta$, we may assume $D$ is effective and $M$ is unchanged. Similarly $\Delta$ must be effective, otherwise $-\Delta$ is effective.  However, $D \cdot \Delta = b \ge 2$ then implies the intersection of the two effective divisors $D$ and $-\Delta$ is negative, so we conclude $D = D' + n(-\Delta)$ for some $D'$ that does not contain $-\Delta$ and some $n >0$. As $D \cdot (-\Delta) <0 $, we must have $(-\Delta)^2 = \Delta^2<0$, contradicting the assumption that $c \ge 0$.  Therefore, we may assume both $D$ and $\Delta$ are effective.  

    Next, we claim that by performing reflections as in Remark \ref{rem:reflections-nef}, we can make $D+\Delta$ nef while preserving the effectivity of $D$ and $\Delta$.  If $R$ is a rational curve such that $(D+ \Delta) \cdot R < 0$, then consider $D' + \Delta'$ where $D' = D + (D \cdot R) R$ and $\Delta' = \Delta + (\Delta \cdot R) R$. By the same remark, $(D')^{2} = D^{2}$ and $(\Delta')^{2} = \Delta^{2}$. Moreover, one can check that
    \[ D' \cdot \Delta' = (D + (D \cdot R) R) (\Delta + (\Delta\cdot R) R) = D \cdot \Delta, \]
    so the intersection matrix of $D'$ and $\Delta'$ is also $M$, and $(D' + \Delta')^2 = (D + \Delta)^2$.  Furthermore, we claim that $D'$ and $\Delta'$ are effective.  We prove this for $D'$ as the proof for $\Delta'$ is identical.
    
    If $D \cdot R \ge 0$, then by definition $D' = D + (D\cdot R) R$ is effective. Now suppose $D \cdot R < 0$ and write $D \cdot R = -k$ where $k>0$. The ideal sheaf sequence for $R$ twisted by $\Oc_X(D)$ implies that $h^0(X,D-R) = h^0(X, D)$ so $D-R$ is effective. By induction, let us assume that $D-nR$ is effective for some $0 < n < k$; we would like to show that $D - (n+1)R$ is also effective. Because $D^2 \ge 0$, we have $D - nR$ is nontrivial for any $n > 0$ (otherwise $D^2 = n^2R^2 = -2n^2$, a contradiction), so Riemann--Roch implies $h^0(X,D-nR)\ge 2 + \frac{1}{2}(D^2 +2kn - 2n^2)$. Consider the sequence
    \[ 0 \to \Oc_X(D - (n+1)R) \to \Oc_X(D  - nR) \to \Oc_R( D- nR ) \to 0 \]
    and note that $\Oc_R(D - nR) = \Oc_{\P^1}(2n-k)$.  If $2n < k$, we have $h^0(X, D-(n+1)R) = h^0(X, D-nR)$ and therefore $D - (n+1)R$ is effective.  For $2n \ge k$, we have $h^0(\P^1, \Oc_{\P^1}(2n-k)) = 2n-k+1$ so
    \begin{align*}
        h^0(X, D - (n+1)R) &\ge 2 + \tfrac{1}{2}(D^2 +2kn - 2n^2) - (2n-k+1) \\
        &= 2 + \tfrac{1}{2}D^2 + n(k-n-2) + k - 1.
    \end{align*}
    If $k \ge n+2$, then we conclude $h^0(X, D - (n+1)R) \ge 2 + \frac{1}{2}D^2$ and hence $D - (n+1)R$ is effective. As we assumed $k > n$, the only remaining case is $k = n+1$, but in this case $n(k-n-2) + k -1 = 0$ so again $h^0(X, D - (n+1)R) \ge 2 + \frac{1}{2}D^2$ and $D - (n+1)R$ is effective. By induction, this proves $D' = D+ (D \cdot R)R = D - kR$ is effective. 

    This demonstrates that $D'$ is effective, and the same argument works for $\Delta'$. Also $D'$ is not trivial because $D' \cdot \Delta' = D \cdot \Delta \ne 0$, so $-D'$ has no sections.  By Riemann--Roch, we conclude that $D'$ has at least 2 sections (and the same holds for $\Delta'$). In summary, after performing finitely many reflections of $D + \Delta$, we may assume $D + \Delta$ is nef and $D, \Delta$ are each effective with at least two sections.  Furthermore, $(D + \Delta)^2 = a + 2b + c > 0$ by assumption so $D + \Delta$ is big.  

    If $D + \Delta $ is not basepoint-free, then by Proposition \ref{prop:base-locus-of-big-and-nef}, $D + \Delta = mE + R$ so case (1) holds.  If $D + \Delta$ is basepoint-free, the general curve $C \in |D + \Delta|$ is smooth and irreducible and we claim that we are in case (2).  Consider the exact sequence 
    \[ 0 \to \Oc_X(-\Delta) \to \Oc_X(D) \to \Oc_C(D) \to 0.\]
    Because $\Delta$ is effective, $h^0(X, - \Delta) = 0$ so we have $h^0(C, D|_C) \ge h^0(X,D) \ge 2 + \frac{1}{2}D^2 = 2 + \frac{a}{2}$.  Because $D$ is effective, we have $h^2(X,D) = h^0(X,-D) = 0$ and therefore $h^1(C, D|_C) \ge h^2(X,-\Delta) = h^0(X, \Delta) \ge 2$.  This implies that $D|_C$ contributes to $\Cliff(C)$ and therefore \[\Cliff(C) \le \Cliff(D|_C) = D \cdot C - 2 (h^0(C, D|_C) - 1) \le a + b - 2(\tfrac{a}{2}+1) = b -2. \]
    If $C$ is not exceptional then $\gon(C) \le b$, and otherwise we conclude $\gon(C) \le b+1$.  
\end{proof}

\subsection{Exceptional curves}

Exceptional curves on K3 surfaces have been classified by Knutsen.  These curves fall into two categories, one specific example constructed by Donagi--Morrison \cite{DonagiMorrison} and an infinite family generalizing the examples of Eisenbud--Lange--Martens--Schreyer \cite{ELMS}. We briefly summarize them here.

\begin{example}\label{ex:DM}
A \textit{Donagi--Morrison example} (cf.\ \cite[(2.2)]{DonagiMorrison})
is a K3 surface presented as the minimal resolution of a double cover branched along a plane sextic with ADE singularities. The composition is denoted $\pi: X \rightarrow \P^2$ and $X$ carries a genus 2 quasi-polarization $L = \pi^{\ast}\Oc(1)$. The arithmetic genus of any curve $C \in |\pi^{\ast}\Oc(3)|$ is 10, and the smooth curves in the codimension $1$ subseries
\[ \pi^* H^0(\P^2, \Oc(3)) \subset H^0(X, \pi^{\ast}\Oc(3)) \]
double cover a plane cubic and hence are bielliptic with gonality $4$. However, a generic member of $|3L|$ maps isomorphically onto its image in $\P^2$, which is a smooth plane sextic, and hence has gonality $5$. Knutsen \cite{Knutsen09} (generalizing work of Ciliberto--Pareschi \cite[Theorem~A]{CP95}) proves that this is the unique situation in which the gonalities of smooth curves in a basepoint-free series on a K3 surface are not constant. 
\end{example}

\begin{remark}
    In the previous example, the branch curve of the map $\pi \colon X \rightarrow \P^{2}$ is allowed to be singular, in contrast with the usual definition of the Donagi--Morrison example considered in \cite{Knutsen09}. The proof of \cite[Theorem 1]{Knutsen09} asserts that a pair $(S,L)$, where $L = 3B$ with $B^2 = 2$, such that the gonalities of smooth curves in $|L|$ is not constant must be the Donagi--Morrison example in \cite{Knutsen09}. Given such an $(S,L)$, it is a priori possible that $B$ fails to be basepoint-free, in which case \Cref{prop:base-locus-of-big-and-nef} implies that $B = 2E + R$ and the curves in the linear series $L$ are degree 3 multisections of the elliptic fibration coming from $|E|$ on $S$. However, the smooth curves must all have gonality 3 (one can see this in multiple ways, e.g. the Castelnuovo--Severi inequality) and hence the gonality is constant in the linear series, which is a contradiction. Thus, $|B|$ is basepoint-free, but not necessarily ample, and there is no reason the branch curve of $\pi$ must be smooth. Therefore, the conclusion of \cite[Theorem 1]{Knutsen09} should allow for the more general Donagi--Morrison examples considered above. 
\end{remark}

The theorem below will motivate the next example:

\begin{theorem}[{\cite[Proposition~4.1]{Knutsen09}}]
Let $L$ be a line bundle on a K3 surface $X$ such that $L \sim 2 D + \Gamma$ with $D$ and $\Gamma$ smooth curves satisfying $D^2 \ge 2$, $\Gamma^2 = -2$, and $D \cdot \Gamma = 1$. Assume that $X$ does not carry any line bundle $B$ with $0 \le B^2 < D^2$ and $0 < B \cdot L  - B^2 \le D^2$.

Then $|L|$ is basepoint-free and all smooth $C \in |L|$ are exceptional curves of genus $g = 2 D^2 + 2 \ge 6$, Clifford index $\Cliff(C) = D^2 - 1 = \frac{g - 4}{2}$, and Clifford dimension $r = \tfrac{1}{2} D^2 + 1$. Moreover, the Clifford index of $C$ is computed by $\Oc_C(D)$ and by no other line bundle on $C$. 
\end{theorem}

\begin{example}\label{ex:ELMS}
A \textit{generalized ELMS example} is a pair $(X, L)$ that satisfies the assumptions in the theorem above. For a generalized ELMS example $(X,L)$, observe that the lattice $\langle D, D + \Gamma \rangle$ is a Donagi--Morrison matrix as defined in Definition \ref{def:DMlattice}.
\end{example}

As mentioned at the beginning of this subsection, Knutsen showed that the two examples above exhaust all possibilities for exceptional curves that can exist on K3 surfaces.

\begin{theorem}[{\cite[Theorem~2]{Knutsen09}}] \label{thm:knutsen}
Let $C$ be a smooth exceptional curve on a K3 surface $X$. Then $C$ is either a smooth plane sextic belonging to the Donagi--Morrison example or $(X, \Oc(C))$ is a generalized ELMS example.
\end{theorem}

\section{Surfaces with degree of irrationality 2 and hyperelliptic curves}\label{section:Irr2HyperellipticCurves}

The goal of this section is to produce smooth hyperelliptic curves on a surface $S$ with degree of irrationality 2. We begin by producing infinitely many smooth rational curves on a mildly singular rational surface.

\begin{lemma}\label{lemma:infiniteRatCurves}
Let $Y$ be a rational surface with finitely many $A_{1}$-singularities, let $B \subset Y$ be a reduced curve, and let $Z \subset Y$ be a finite set. Then $Y$ contains infinitely many smooth rational curves which intersect $B$ transversely and avoid $Z$.
\end{lemma}

\begin{proof}
We run an MMP on $Y$, contracting $K_{Y}$-negative curves to get a sequence of divisorial contractions
\[ Y = Y_0 \xrightarrow{\psi_0} Y_1 \xrightarrow{\psi_1} Y_2 \to \cdots \to Y_n \]
such that $Y_n$ is either minimal, i.e. $K_{Y_n}$ is nef, $Y_{n}$ Fano with Picard rank 1, or $Y_{n}$ admits a map to $\P^1$ with generic fiber $\P^1$.

We claim that at each step, $Y_{i}$ has at worst finitely many $A_{1}$-singularities. To see this, let us induct on $i$ (with the base case $i = 0$ satisfied by assumption). Suppose that $Y_{i}$ satisfies the inductive hypothesis. Then $\psi_{i}$ corresponds to the contraction of an integral divisor $R$, and Artin's contractibility criterion implies $R^{2} < 0$. Since $Y_{i}$ has at worst $A_{1}$ singularities, it is Gorenstein and hence $K_{Y_{i}} \cdot R \in\Zb$.  By assumption $K_{Y_{i}} \cdot R < 0$ and therefore $K_{Y_{i}} \cdot R \le -1$.  It follows that $R$ is a smooth curve with $K_{Y_{i}} \cdot R = -1$ (see, e.g. \cite[Lemma 4.8]{DeVleming}). The local class group at each node is $\Zb/2$, so $2R$ is Cartier and hence $2R \cdot R \in \Zb$, which implies that $R^{2} \in \frac{1}{2} \Zb$.  Using the adjunction formula
\[ 2 p_{a}(R) - 2 + \sum_{\text{sing pts on } R} (1 - \tfrac{1}{2}) = (K_{Y_i}+R) \cdot R, \]  we conclude \[-2 + \tfrac{n}{2} = K_{Y_{i}} \cdot R + R^{2} \leq -\tfrac{3}{2}, \]
where $n$ is the number of $A_1$ singularities on $R$.  It follows that either $n = 0$ and $R^{2} = -1$, or $n = 1$ and $R^{2} = -\frac{1}{2}$. It is straightforward to show that, under either operation, the image of $R$ is a smooth point. 
Observe that the number of singular points of $Y_{i+1}$ is then at most the number of singular points of $Y_{i}$, and $Y_{i+1}$ has only $A_1$ singularities. 

We conclude $Y_{n}$ is a Gorenstein rational surface with finitely many $A_{1}$-singularities. In particular, $K_{Y_{n}}$ is not pseudoeffective so $K_{Y_n}$ is not nef. If $Y_{n}$ admits a fibration over $\P^1$ with generic fiber $\P^1$, then a general fiber of the fibration $Y_{n} \rightarrow \P^{1}$ is a smooth rational curve whose intersection with the image of $B$ is transverse (if nonempty) and misses any finite set $Z$, so we are done. If $Y_n$ is Fano of Picard rank 1, the classification of Del Pezzo surfaces of rank 1 (see, for example \cite{Lacini}) yields the only possibilities for which $Y_{n}$ with at worst $A_1$ singularities are either $\P^{2}$ or $\P(1, 1, 2)$. For $\P^{2}$, one may take a general line, which will meet any fixed curve $B$ transversely and avoid any finite set. For $\P(1,1,2)$, which is the quadric cone, a general hyperplane section gives a smooth conic, which will meet $B$ transversely and avoid any finite set.  The proof is complete.
\end{proof}

\begin{corollary}
    Let $X$ be a smooth surface with $\irr(X) = 2$. Then $X$ contains infinitely many smooth curves which are either elliptic or hyperelliptic.
\end{corollary}

\begin{proof}
    Suppose $\kappa(X) = -\infty$. Then $X \cong_{\text{bir.}} C \times \P^{1}$ for some curve $C$ with $\gon(C) = 2$. By Tsen's theorem, the Albanese morphism $\alpha \colon X \rightarrow C$ has a section. But then the generic fiber of $\alpha$ is isomorphic to $\P^{1}_{\C(C)}$, so there are infinitely many sections of $\alpha$, which are smooth curves isomorphic to $C$. Now suppose $\kappa(X) \geq 0$, and let $X'$ be the minimal model of $X$. Then $X'$ is smooth and $K_{X'}$ is nef, so $\Bir(X') = \Aut(X')$ by the uniqueness of minimal models. Given a degree 2 rational map $X \dashrightarrow \P^{2}$, we get an involution $\phi \in \Bir(X)$ which descends to an involution $\phi' \in \Bir(X') = \Aut(X')$. The quotient $Y \colonequals X'/ \langle \phi' \rangle$ is thus a normal rational surface with at worst $A_{1}$-singularities. Let $B$ be the branch divisor of $Y$ and let $Z \subset Y$ be the finite set consisting of singular points of $Y$ plus the images of $(-1)$-curves under the composition $\phi \colon X \rightarrow X' \rightarrow Y$. By Lemma \ref{lemma:infiniteRatCurves}, there exist infinitely many smooth rational curves in $Y$ which are transverse to $B$ and avoid $Z$. Their pre-images under $\phi$ are therefore smooth curves of gonality equal to $2$.
\end{proof}

We will use this to give equivalent conditions for a K3 surface $X$ to have $\irr(X) = 2$. First, let us recall the following definition.

\begin{definition}[Saint-Donat, Reid]\label{def:hyperelliptic}
Let $L \in \Pic{X}$ be a line bundle on a K3 surface $X$ such that $L^{2} > 0$ and $|L|$ has no fixed components. If the associated map $X \rightarrow \P^{g}$ is 2-to-1 onto its image (which is a surface of degree $g-1$), then $|L|$ is called a \textit{hyperelliptic linear series}.
\end{definition}

The following theorem is closely related to the Enriques--Campedelli Theorem \cite{Dolgachev73, Reid76a} and we believe it is well-known to experts, but we could not find a suitable reference. Thus, we will include a short proof.

\begin{theorem}\label{thm:irr2hypcurve}
Let $X$ be a smooth projective K3 surface. The following are equivalent:
\begin{enumerate}[label=(\alph*)]
\item $\irr(X) = 2$;
\item there exists a basepoint-free line bundle $L$ such that $|L|$ is hyperelliptic;
\item $X$ contains a smooth hyperelliptic curve of genus $g \geq 2$;
\item $X$ contains a smooth curve of genus $g = 2$ or there is a genus 1 fibration $X \rightarrow \P^{1}$ with an irreducible multisection of degree $\leq 2$.
\end{enumerate}
\end{theorem}

\begin{proof}

$(b) \implies (a)$: The image of $X$ under $|L|$ as in Definition \ref{def:hyperelliptic} is a surface $S \subset \P^{g}$ of minimal degree which is necessarily rational.

$(c) \implies (b)$: Let $C$ be the smooth curve of genus $g \geq 2$. Then by adjunction, the restriction of the linear series $|C|$ to any smooth curve $C' \in |C|$ is the canonical linear series for $C'$. The claim now follows from the well-known fact that $C$ is hyperelliptic if and only if every smooth member $C' \in |C|$ is hyperelliptic. A proof of this fact can be found in \cite[Theorem 3.2]{DonagiMorrison}, who attributed it to \cite{Saint-Donat74}.

$(d) \implies (a)$: the first part is trivial. The second part essentially follows from \cite[Proposition 1]{Yoshihara}, which we will summarize. Suppose $X$ admits a genus 1 fibration $\pi \colon X \rightarrow \P^{1}$ with an irreducible multisection $\Gamma$ of degree $\leq 2$, and let $E$ be the general fiber of $\pi$. By doubling $\Gamma$ if necessary, we may assume that $\Gamma \cdot E = 2$. Then $\pi_*\Oc_X(\Gamma)$ is a rank two vector bundle on $\P^1$ and there is an induced rational map $X \dashrightarrow \P(\pi_*\Oc_X(\Gamma))$. By construction, this map has degree 2 and $\P(\pi_*\Oc_X(\Gamma))$ is a rational surface, so we are done. 

$(a) \implies (c)$ and $(d)$ will be given simultaneously. Given a degree 2 map $X \dashrightarrow \P^{2}$, let $\sigma \in \Bir(X) = \Aut(X)$ be the induced involution. Since $\sigma$ is anti-symplectic \cite{Nikulin81}, the quotient $Y \colonequals X/\sigma$ is a smooth rational surface with smooth branch locus $B \in |-2K_{Y}|$. By Lemma~\ref{lemma:infiniteRatCurves}, there are infinitely many smooth rational curves in $Y$ whose preimages in $X$ are smooth curves of gonality equal to 2. Moreover, in the proof of the Lemma $Y^{\min}$ is either isomorphic to $\P^{2}$ or a Hirzebruch surface $\F_{n}$ (with $0 \leq n \leq 4$ by \cite{Reid76a}).

When $Y^{\min} = \P^{2}$, the preimages under $X \rightarrow Y^{\min}$ of general lines or conics become smooth curves of genus 2 or smooth curves of genus 5, respectively. For $Y^{\min} = \F_{n}$, the genus 1 fibration on $X$ comes from pulling back a general fiber class of $\F_{n}$ to get a basepoint-free pencil on $X$, and the multisection of degree $\leq 2$ arises from taking the reduced preimage of any section. On the other hand, we also claim that in this case $X$ contains smooth hyperelliptic curves of genus $\geq 2$. To see this, let $S \subset \F_{n}$ denote the unique section of self-intersection $-n$ if $n > 0$ or any section of self-intersection $0$ if $n = 0$, and let $F$ be a ruling. Then for $m \ge \max\{1,n\}$, the linear series $|S + mF|$ is basepoint-free and a general element meets the branch locus of $X \rightarrow Y^{\min}$ transversely in at least $6$ points. As these are sections of the ruled surface, their preimages on $X$ are smooth hyperelliptic curves of genus $g \ge 2$.
\end{proof}

\begin{example}\label{ex:irrPicZK3}
Let $X$ be a projective K3 surface with $\Pic(X) = \Zb \cdot L$. By \cite[Lemma 5.2]{EO15}, it follows that 
\[ \Bir(X) = \Aut(X) = \begin{cases} \{ \text{id}_{X} \} & \text{if } L^{2} \not= 2 \\ \Zb/2\Zb & \text{if } L^{2} =  2. \end{cases} \]
Thus, $\irr(X) > 2$ when $L^{2} \not= 2$.
\end{example}

\section{Lattice theoretic characterizations of K3 surfaces with $\irr = 2$}\label{section:hyperellipticK3s}

The main goal of this section is to prove the following theorem.

\begin{theorem}\label{thm:Irr2K3s}
Let $X$ be a projective K3 surface. Then $\irr(X) = 2$ if and only if one of the following conditions holds:
\begin{enumerate}[label=(\alph*)]
    \item there exists an ample class $L \in \Pic{X}$ with $L^{2} = 2$;
    \item One of the lattices embeds primitively into $\Pic{X}$:
    \[ \begin{pmatrix} 0 & 2 \\ 2 & -2 \end{pmatrix}, \quad \begin{pmatrix} 0 & 2 \\ 2 & 0 \end{pmatrix}, \quad \text{or} \quad \begin{pmatrix} 0 & 1 \\ 1 & -2 \end{pmatrix}. \]
\end{enumerate}
\end{theorem}

\noindent Note that these conditions are not mutually exclusive. Up to a change of basis, the third lattice is isomorphic to the hyperbolic plane.

\begin{remark}\label{rmk:which-lattice}
    By the analysis in this section, we will see that if one of the lattices above embeds into $X$, then in fact $X$ can always be realized as a double cover of the following surfaces:
    \begin{enumerate}
        \item $\begin{pmatrix} 0 & 2 \\ 2 & 0 \end{pmatrix}$ implies that K3 can be constructed as a double cover of $\P^1 \times \P^1$;
        \item $\begin{pmatrix} 0 & 2 \\ 2 & -2 \end{pmatrix}$ implies that $X$ can be constructed as a double cover of $\F_1$; 
        \item $\begin{pmatrix} 0 & 1 \\ 1 & -2 \end{pmatrix}$ gives an elliptic K3 with a section, which can always be realized as a double cover of $\F_4$. 
    \end{enumerate}
\end{remark}

\subsection{Existence of lattice embedding implies degree of irrationality 2}

In this subsection, we prove the reverse direction of Theorem \ref{thm:Irr2K3s}.

\begin{lemma}\label{lem:K3DoubleCoverP2}
Let $X$ be a K3 surface over an algebraically closed field of characteristic $\not= 2$. Then there exists an involution $\iota$ on $X$ with quotient $X / \iota \cong \P^{2}$ if and only if there exists an ample class $L \in \Pic{X}$ with $L^{2} = 2$.
\end{lemma}

\begin{proof}
The forward direction is obvious since one can choose $L$ to be the pullback of the hyperplane class on $\P^{2}$. For the reverse direction, the condition that $L$ is ample implies by \cite[Proposition 8.1]{Saint-Donat74} that either $L \cong \Oc_{X}(C)$ for some irreducible curve $C$ with $p_{a}(C) > 1$, or $L \cong \Oc_{X}(k \cdot E + \Gamma)$ with $k \geq 3$, where $E$ is an irreducible curve with $p_{a}(E) = 1$ and $\Gamma$ is a smooth rational $(-2)$-curve with $\Gamma \cdot E = 1$.  The latter is impossible since we are assuming $L^{2} = 2$. So $L \cong \Oc_{X}(C)$ is basepoint-free \cite[Chapter 2, Proposition 3.5]{Huybrechts16} and defines a morphism $\pi \colon X \rightarrow \P^{2}$. $L$ being ample tells us that $\pi$ is finite of degree 2.  This gives the desired involution $\iota$.
\end{proof}

\begin{lemma}\label{lem:rank2irr2}
Let $X$ be a projective K3 surface, and suppose there exists a primitive embedding of one of the following lattices into the Picard lattice $\Pic{X}$:
\[ \begin{pmatrix} 0 & 2 \\ 2 & -2 \end{pmatrix}, \quad \begin{pmatrix} 0 & 2 \\ 2 & 0 \end{pmatrix}, \quad \text{or} \quad \begin{pmatrix} 0 & 1 \\ 1 & -2 \end{pmatrix}. \]
Then $\irr(X) = 2$.
\end{lemma}

\begin{proof}

We will first prove the theorem under the assumption that $\Pic{X}$ is equal to one of the three lattices listed above. Let $A$ and $B$ denote the standard generators for each of these lattices. For example, in each lattice above we have $A^{2} = 0$. Recall from Riemann--Roch that
\[ h^{0}(X, A) + h^{0}(X, -A) = h^{0}(X, A) + h^{2}(X, A) \geq \frac{1}{2} A^{2} + 2 = 2. \]
Thus, up to replacing $A$ with $-A$, we may assume that $A$ is effective. Replacing $B$ with $-B$ if necessary, we may also assume that $B$ is effective. However, a priori $A \cdot B$ could flip signs.

Suppose that $X$ does not contain any $(-2)$-curves. Then $A \cdot B \geq 0$ so we have the same lattices above, and $|A|$ does not contain any fixed components so it is basepoint-free and defines a pencil $\pi \colon X \rightarrow \P^{1}$. For $m$ large enough, the class $B + mA$ is effective, basepoint-free, and defines a multisection of degree at most 2. By \Cref{thm:irr2hypcurve}, we are done.

Now suppose $X$ contains a $(-2)$-curve; then by the rank assumption on $\Pic X$, there is at most one such curve. We can exclude the second case $\begin{pmatrix} 0 & 2 \\ 2 & 0 \end{pmatrix}$ because for any $a, b \in \Zb$ we arrive at a contradiction when setting $-2 = (aA + bB)^{2} = 4ab$. In the first case where $\begin{pmatrix} 0 & 2 \\ 2 & -2 \end{pmatrix}$ is the intersection matrix, let us write $aA + bB$ for the $(-2)$-curve. Then
\[ -2 = (aA + bB)^{2} = 4ab + b^{2} (-2) \implies -1 = (2a - b)b. \]
So $b = \pm 1$ and $a = 0$, which implies that either $B$ or $-B$ is a $(-2)$-curve. After choosing the correct sign, we may assume that $B$ is the class of a smooth $(-2)$-curve. In this process, it may be possible for the intersection $A \cdot B$ to be $-2$. If $(A \cdot B) < 0$, then one can apply a reflection to replace $A$ with $A + (A \cdot B) B$. Note that
\[ (A + (A \cdot B) B) \cdot B = A \cdot B + (A \cdot B) (-2) = -(A \cdot B), \]
so $A$ is nef and we arrive back at the same intersection matrix. For a nef class $A$ with $A^{2} = 0$, \cite[Chapter 2, Proposition 3.10]{Huybrechts16} implies that $A$ is basepoint-free and of the form $k \cdot A'$ where $A'$ is the class of a smooth elliptic curve on $X$. The condition $A \cdot B = 2$ implies that $k = 1$ or $2$, and so we are done by \Cref{thm:irr2hypcurve}. In the third case where $\begin{pmatrix} 0 & 1 \\ 1 & -2 \end{pmatrix}$ is the intersection matrix,
\[ -2 = (aA + bB)^{2} = 2ab \cdot (A \cdot B) + b^{2}B^{2} = 2b (a - b). \]
The only solutions are $\{ b = \pm 1, a = 0 \}$. By a similar argument as above, there exists an elliptic pencil with a degree two multisection so by \Cref{thm:irr2hypcurve} we have $\irr(X) = 2$.

Suppose $X$ is now a K3 surface and there exists a primitive embedding $M \hookrightarrow \Pic X$, where $M$ is one of the three lattices in the hypothesis. There is a sub-locus $\Omega_{M}$ of the period domain consisting of $[\omega] \in \P(\Lambda_{\text{K3}} \otimes \C)$ such that every class in $M$ is algebraic, which is equivalent to the fact that the pairing $\langle \omega, \alpha \rangle = 0$ for all $\alpha \in M$. The very general element in $\Omega_{M}$ corresponds (via the local Torelli theorem) to a K3 surface with Picard lattice equal to $M$, since containing an additional algebraic class is a linear condition in $\Omega_{M}$. By surjectivity of the period map, there is a small deformation $X'$ of $X$ whose Picard group is precisely $\Pic{X'} = M$. But then \cite[Proposition C]{CSirr} implies that $2 = \irr(X') \geq \irr(X)$ and so $\irr(X) = 2$.
\end{proof}

\subsection{From double covers to lattices}

Now we will prove the other direction of Theorem \ref{thm:Irr2K3s}. Most of this relies on the study of curves on Hirzebruch surfaces, as K3 surfaces with $\irr=2$ are (up to taking a resolution) double covers of $\P^2$ or a Hirzebruch surface $\F_n$ for $0 \le n \le 4$ by \cite{Reid76a}. The proof of the following lemma will be omitted.

\begin{lemma}\label{lem:dimension-count}
    The moduli space of pairs $(\F_n,B)$ where $B \in |-2K_{\F_n}|$ is a curve with at worst ADE singularities (up to equivalence by automorphism of $\F_n$) has dimension $18$ if $n = 0,1,4$, has dimension $17$ if $n = 2,3$, and is empty if $n \ge 5$. Furthermore, the pairs $(\F_2, B)$ and $(\F_3,B)$ are, respectively, degenerations of the pairs $(\F_0, B)$ and $(\F_1,B)$.
\end{lemma}

\begin{lemma}\label{lem:doublecover-lattice}
    Let $Y$ be the double cover of a Hirzebruch surface $\F_n$, $0 \le n \le 4$, branched over a divisor $B \in |-2K_{\F_n}|$ such that $B$ has at worst ADE singularities, and denote by $\phi_Y : Y \to \F_n$.  Let $X$ be the K3 surface that is the minimal resolution of $Y$, and denote by $\phi_X: X \to \F_n$ the generically finite map. Then, there exists a primitive embedding of one of the following lattices into the Picard lattice $\Pic{X}$:
\[ \begin{pmatrix} 0 & 2 \\ 2 & -2 \end{pmatrix}, \quad \begin{pmatrix} 0 & 2 \\ 2 & 0 \end{pmatrix}, \quad \text{or} \quad \begin{pmatrix} 0 & 1 \\ 1 & -2 \end{pmatrix}. \]
\end{lemma}

\begin{proof}
Let $\pi: \F_n \to \P^1$ be the ruled surface map and let $f \subset \F_n$ be a general fiber of $\pi$.  For $n \ge 1$, let $s \subset \F_n$ be the negative section such that $s^2 = -n$, and for $n = 0$ let $s \subset \F_0 = \P^1 \times \P^1$ be a general horizontal section of $\pi$ with $s^2 = 0$.  

Let $E \subset X$ be the pullback of the general fiber $f$ and let $D \subset X$ be a reduced component of the pullback of $s$ mapping finitely to $s$.  We claim that the intersection matrix of $E$ and $D$ is one of \[ \begin{pmatrix} 0 & 2 \\ 2 & -2 \end{pmatrix}, \quad \begin{pmatrix} 0 & 2 \\ 2 & 0 \end{pmatrix}, \quad \text{or} \quad \begin{pmatrix} 0 & 1 \\ 1 & -2 \end{pmatrix}. \]
Because $E$ is the pullback of a general fiber, we have $E^2 = 0$. To determine the intersection matrix of $E$ and $D$, we must distinguish the cases when $s$ is either contained in the branch locus of $\phi_X$ or not.  For $n = 0$, as $s$ was a general horizontal section, it is not contained in the branch locus but meets the branch locus in four smooth points.  The pullback is a unique curve $D$ mapping with degree 2 to $s$, and $s^2 = 0$ implies that $D^2 = 0$, and as $s \cdot f = 1$, the projection formula yields $E \cdot D = 2$.  Therefore, the intersection matrix of $E$ and $D$ is $\begin{pmatrix} 0 & 2 \\ 2 & 0 \end{pmatrix}$. 

Now assume $n \ge 1$ and recall that $-2K_{\F_n} \equiv 4s + (2n+4)f$.  Assume first that $s$ is not contained in branch locus $B$ of $\phi_X$, which necessarily implies that $n = 1,2$ by the analysis of the linear series $|-2K_{\F_n}|$.  Then, $B \cdot s = 4 - 2n$, so when $n = 1$, $B \cdot s = 2$ and when $n = 2$, $B \cdot s = 0$.  When $n = 2$, the pullback of $s$ is the union of two disjoint $(-2)$-curves.  Let $D$ be one such curve.  Then, $D^2 = -2$ and $E \cdot D = 1$, and therefore the intersection matrix of $E$ and $D$ is  $\begin{pmatrix} 0 & 1 \\ 1 & -2 \end{pmatrix}$. When $n =1$ and $s$ meets the branch locus in two smooth points, the pullback of $s$ is irreducible and equal to $D$, and therefore by the projection formula we have $D^2 = -2$ and $E \cdot D = 2$, so the intersection matrix of $E$ and $D$ is $\begin{pmatrix} 0 & 2 \\ 2 & -2 \end{pmatrix}$.  Next we consider the case when $n = 1$ and $s$ meets the branch locus either tangentially at a smooth point or at a singular point. In the first case, the pullback of $s$ consists of two $(-2)$-curves meeting at a single point, and taking $D$ to be either curve, from the projection formula, we get the intersection matrix of $E$ and $D$ to be $\begin{pmatrix} 0 & 1 \\ 1 & -2 \end{pmatrix}$.  In the second case, by Lemma \ref{lem:factors-through-min-res}, $X$ is the minimal resolution of a double cover of the blow-up of the singular point of $B$ on $s$, branched over the strict transform of $B$.  In this case, the pullback of $s$ to $X$ is the union of two disjoint $(-2)$-curves and taking either one for $D$, we get the intersection matrix $\begin{pmatrix} 0 & 1 \\ 1 & -2 \end{pmatrix}$.
    
Next, assume $n \ge 1$ and $s$ is contained in the branch locus $B$ of $\phi_X$.  Write $B = s + B'$, where $B' \equiv 3s + (2n+4)f$ is reduced, has at worst $ADE$ singularities and does not contain $s$.  Then, $B' \cdot s  = 4 - n$.  If $s$ meets $B'$ transversely at $4-n$ distinct points, the pullback of $s$ on $Y$ is a doubled irreducible curve passing through $4-n$ $A_1$ singularities.  Letting $D_Y$ be the reduced preimage of $s$ on $Y$, by the projection formula $D_Y^2 = -\frac{n}{2}$ and therefore its strict transform $D \subset X$ in the minimal resolution of $Y$ has $D^2 = -2$.  From the projection formula $E \cdot \phi_X^*(s)  = 2$, and as $\phi_X^*(s)$ is doubled along $D$, we have $E \cdot D = 1$, so the intersection matrix of $E$ and $D$ is  $\begin{pmatrix} 0 & 1 \\ 1 & -2 \end{pmatrix}$.  

    It remains to consider the case when $n =1,2$ and $s$ does not meet $B'$ transversely.  When $n = 2$, either $s$ meets $B'$ tangentially at a smooth point of the curve $B'$ or $s$ meets $B'$ at a singular point of $B'$.  We apply Lemma \ref{lem:factors-through-min-res} in either case to say that there is a generically finite map $X \to S'$ where $S'$ is a partial log resolution of the pair $(\F_n, B)$.  In the first case, blowing up $s$ twice resolves the intersection of $s$ and $B'$, and the exceptional divisors of these blow ups are not part of the branch locus on the new surface, denoted $S'$.  The strict transform of $s$ on $S'$ is then a $-4$ curve that is contained in the branch locus of the map $X \to S'$, but meets no other component of the branch locus.  Therefore, its reduced preimage on $X$ is a $(-2)$-curve $D$ that meets the preimage $E$ of a general fiber in one point, and the corresponding intersection matrix of $E$ and $D$ is $\begin{pmatrix} 0 & 1 \\ 1 & -2 \end{pmatrix}$.  In the second case, $s$ meets $B'$ at a singular point (to order 2) and the intersection of $s$ and $B'$ is resolved after one blow up $S' \to \F_2$, and the corresponding exceptional divisor of the blow up is contained in the branch locus of $X \to S'$.  Blowing up the intersection of the exceptional divisor and the strict transform of $s$ to get a new surface $S'' \to S'$, the second exceptional divisor is not contained in the branch locus of $X \to S''$.  Therefore, the strict transform of $s$ on $S''$ is a $-4$ curve that does not meet any other component of the branch locus, and hence its reduced preimage on $X$ is a $(-2)$-curve $D$ that meets $E$ in one point, and the corresponding intersection matrix of $E$ and $D$ is $\begin{pmatrix} 0 & 1 \\ 1 & -2 \end{pmatrix}$.

    Finally, we treat the case of $n = 1$ again using Lemma \ref{lem:factors-through-min-res}.  In this case, either $s$ meets $B'$ transversely at one point and tangentially at a second smooth point of $B'$, or $s$ meets $B'$ transversely at one point and at a singular point of $B'$ to order 2, or $s$ meets $B'$ at a single point to order 3.  By the same argument as above, in the first or second case, blowing up the transverse intersection once and the order 2 intersection twice results in $X \to S'$ where the strict transform of $s$ is a $-4$ curve contained in the branch locus that does not meet any other component of the branch locus.  Its reduced preimage on $X$ is a $(-2)$-curve $D$ that meets $E$ in one point, and the corresponding intersection matrix of $E$ and $D$ is $\begin{pmatrix} 0 & 1 \\ 1 & -2 \end{pmatrix}$.  In the final case when $s$ meets $B'$ to order 3, after blowing up the intersection point three times, we obtain a surface $S'$ and generically finite map $X \to S'$ where the strict transform of $s$ is a $-4$ curve contained in the branch locus that does not meet any other component of the branch locus and the same argument yields the intersection matrix of $E$ and $D$ on $X$ is $\begin{pmatrix} 0 & 1 \\ 1 & -2 \end{pmatrix}$.

    If the third matrix $\begin{pmatrix} 0 & 1 \\ 1 & -2 \end{pmatrix}$ embeds into $\Pic(X)$, then it must be primitive since it has determinant 1. For the other two matrices $\begin{pmatrix} 0 & 2 \\ 2 & 0 \end{pmatrix}$ and $\begin{pmatrix} 0 & 2 \\ 2 & -2 \end{pmatrix}$, one can show that if the embedding in $\Pic X$ is not primitive, then the saturation is in fact the third lattice. 
\end{proof}

The previous lemma made use of the following, which is justified in \cite[Section 7]{MR1913513}.

\begin{lemma}\label{lem:factors-through-min-res}
    Let $\phi_Y : Y \to S$ be the double cover of a smooth surface $S$ branched over a divisor $B$ such that $B$ has at worst ADE singularities.  Let $X$ be the minimal resolution of $Y$, and denote by $\phi_X: X \to S$ the generically finite map. Let $\widetilde{\psi}: \widetilde{S} \to S$ be the minimal log resolution of $(S,B)$ and define $\widetilde{B}$ by $K_{\widetilde{S}} + \tfrac{1}{2}\widetilde{B} = \widetilde{\psi}^*(K_{S} + \tfrac{1}{2}B)$. Now let $\overline{\psi}: \overline{S} \to \widetilde{S}$ be the blow up of each node of $\widetilde{B}$ and let $\overline{B}$ be the strict transform of $\widetilde{B}$.  Note that $K_{\overline{S}} + \tfrac{1}{2}\overline{B} = \overline{\psi}^*(K_{\widetilde{S}} + \tfrac{1}{2}\widetilde{B})$. Let $\psi: \overline{S} \to S$ be the composition and let $\overline{X}$ be the double cover of $\overline{S}$ branched over $\overline{B}$.
    
    Then $\overline{X} = X$, and given any partial resolution $\overline{S} \to S' \to S$ with associated double cover $X'$, $\phi_X: X \to S$ factors through the generically finite map $X' \to S$.
\end{lemma}

The forward direction of Theorem~\ref{thm:Irr2K3s} now follows from the fact that \Cref{thm:irr2hypcurve}(b) gives the existence of a hyperelliptic line bundle on any K3 surface $X$ with $\irr(X) = 2$. By work of Reid \cite{Reid76a}, $X$ is a double cover of $\P^2$ or $\F_n$ with $0 \le n \le 4$.  Furthermore, if $X$ is a double cover of $\P^2$, we may assume the hyperelliptic line bundle is ample (if it is only big and nef, the cover factors through a map to $\F_1$).  Thus, if $\Pic{X}$ does not contain an ample class of square $2$, we can apply \Cref{lem:doublecover-lattice} to conclude. This completes the proof of Theorem~\ref{thm:Irr2K3s}.

\subsection{An equivalent formulation} While Theorem \ref{thm:Irr2K3s} and Remark~\ref{rmk:which-lattice} give a geometric characterization of the K3 surfaces with $\irr = 2$, there is an equivalent criterion based on the existence of \textit{arbitrary} classes of square 2.

\begin{theorem}\label{thm:Irr2K3s-version2}
    Let $X$ be a projective K3 surface.  Then $\irr(X) = 2$ if and only if $\Pic{X}$ contains a class of square 2 or $\begin{pmatrix} 0 & 2 \\ 2 & 0 \end{pmatrix}$ embeds primitively in $\Pic{X}$. 
\end{theorem}

The proof relies on the following result.

\begin{lemma}\label{lem:class-of-square-2}
    Let $(X,L) \in \Fc_g$ be a quasi-polarized K3 surface of genus $g \ge 3$.  If $\Pic{X}$ contains a class of square 2, then $X$ is contained in a Noether--Lefschetz divisor corresponding to a lattice
    \[ \Lambda = \begin{pmatrix} 2 & e \\ e & 2g-2 \end{pmatrix} \]
    for some $e>0$ such that $e^2 > 4g-4$, and moreover $\irr(X) = 2$.
\end{lemma}

\begin{proof}
    Let $D$ represent the class of square 2.  Up to replacing $D$ with $-D$ and performing reflections as in Remark \ref{rem:reflections-nef}, we may assume $D$ is nef.  As $D^2 = 2$, it is big and necessarily primitive.  Let $\Lambda = \begin{pmatrix} 2g-2 & e \\ e & 2 \end{pmatrix} $ be the intersection matrix of $L$ and $D$, which embeds in $\Pic{X}$.  By the Hodge index theorem and nefness of $D,L$, we conclude $e >0$ and $e^2 > 4g - 4$.  This proves the first statement that $X$ is contained in the given Noether--Lefschetz divisor.   

    For the statement on $\irr(X)$, by \cite[Corollary 3.1]{CE25} it suffices to show that the very general point $(X,L)$ of the Noether--Lefschetz divisor has $\irr(X) = 2$. Because $(X,L)$ is very general, it follows that $\Pic{X} = \Lambda$.  First note that if $D$ is basepoint free, then $D$ defines a degree 2 map to $\P^{2}$ and hence $\irr(X) = 2$.  Supposing $D$ is not basepoint-free, we will show that $\Lambda$ is actually equivalent to either $\begin{pmatrix} 0 & 1 \\ 1 & -2 \end{pmatrix} $ or $\begin{pmatrix} 0 & 2 \\ 2 & -2 \end{pmatrix} $.  In either case, by Theorem \ref{thm:Irr2K3s}, it will follow that $\irr(X) = 2$.  
    
    By the assumption that $D$ is not basepoint-free, by Proposition \ref{prop:base-locus-of-big-and-nef}, $\Pic{X} = \Lambda$ must represent both $0$ and $-2$.  Suppose $E = a L + b D$ satisfies $(a L + b D)^{2} = 0$ for $a,b \in \Zb$. Expanding the expression gives
\[
E^2 = (a L + b D)^{2} = (2g-2) a^{2} + 2 b^{2} + 2abe = 0.
\]
Dividing out by $2a^{2}$ and setting $\lambda \colonequals \frac{b}{a}$ yields $\lambda^{2} + e \lambda + (g-1) = 0$. The discriminant of this quadratic equation is $e^{2} - 4 (g-1) = - \det \Lambda$, so we see that we must have $e^{2} - 4(g-1)$ a perfect square by existence of the rational solution $\lambda$.  Write $e^2  - 4(g-1) = d^2$ for $d > 0$. 

Next, choose a primitive class $E$ with $E^2 = 0$ and extend it to a basis $E,F$ of the lattice $\Lambda$.  Up to switching the sign of $F$, we may assume the intersection matrix of $E,F$ is given by $\Lambda' = \begin{pmatrix} 0 & s \\ s & t \end{pmatrix}$ with $s > 0$, and $s^2 = -\det \Lambda' = -\det \Lambda = d^2$, so $s = d$.  Note that $t$ is even as $\Lambda'$ is even.  

By definition, $\Lambda$ (and hence $\Lambda'$) represents $2$, so we must have a solution to
\[ 2=  (xE + yF)^2 = 2dxy + ty^2  = y(2dx + ty).\]
Because $t$ is even, the second factor is even, so we must have $y = \pm 1$ and $2dx + ty = \pm 2$, hence $t \equiv  2 \pmod{2d}$.  However, by assumption, $\Lambda$ also represents $-2$, so by the same argument, we must have $t \equiv -2 \pmod {2d}$.  This is possible if and only if $d = 1$ or $d =2$

Furthermore, because $(F + kE)^2 = t + 2kd$ and $(F+kE) \cdot E = d$, up to replacing $F$ with $F + kE$ for some $k \in \Zb$ we may assume $t = -2$, i.e. $\Lambda' = \begin{pmatrix} 0 & d & d & -2 \end{pmatrix}$.  Because $d =1$ or $d =2$, we have shown that $\Lambda$ is actually equivalent to either $\begin{pmatrix} 0 & 1 \\ 1 & -2 \end{pmatrix} $ or $\begin{pmatrix} 0 & 2 \\ 2 & -2 \end{pmatrix} $ and hence $\irr(X) = 2$, as desired. 
\end{proof}

With Lemma \ref{lem:class-of-square-2}, it is straightforward to prove Theorem \ref{thm:Irr2K3s-version2}. 

\begin{proof}[Proof of \Cref{thm:Irr2K3s-version2}]
    Assume $\irr(X) = 2$.  By Theorem \ref{thm:Irr2K3s}, it suffices to show that the lattices $\begin{pmatrix} 0 & 1 \\ 1 & -2 \end{pmatrix} $ and $\begin{pmatrix} 0 & 2 \\ 2 & -2 \end{pmatrix} $ contain a class of square 2.  Denoting the classes spanning the lattices by $E, R$, the required class is given by $2E +R$ in the first case and $E+R$ in the second. 

    To complete the proof, assume  $\Pic{X}$ contains a class of square 2 or $\begin{pmatrix} 0 & 2 \\ 2 & 0 \end{pmatrix}$ embeds primitively in $\Pic{X}$.  In the first case, Lemma \ref{lem:class-of-square-2} implies $\irr(X) = 2$, and in the second case, Theorem \ref{thm:Irr2K3s} implies $\irr(X) = 2$.  
\end{proof}

Finally, the following result immediately implies \Cref{thm:Hg-in-G2g}.

\begin{theorem}\label{thm:codim2-hyperelliptic}
    For any fixed even genus $g \ge 4$, the locus $\Hc_g$ is a countably infinite union of divisors. For any fixed odd genus $g \ge 3$, $\Hc_{g}$ is a countably infinite union of divisors and codimension 2 components parametrizing K3s which arise as double covers of $\P^1 \times \P^1$ such that the quasi-polarization is NOT pulled back from $\P^1 \times \P^1$. 
\end{theorem}

\begin{proof}
By Lemma \ref{lem:class-of-square-2}, there exist countably infinitely many Noether--Lefschetz divisors parametrizing K3 surfaces with a class of square 2.  By Theorem \ref{thm:Irr2K3s-version2}, these are all contained in $\Hc_g$. On the other hand, if a K3 surface in $\Hc_{g}$ does not contain a class of square two, then by Theorem \ref{thm:Irr2K3s-version2}, the matrix $\begin{pmatrix} 0 & 2 \\ 2 & 0 \end{pmatrix}$ embeds primitively in $\Pic{X}$.

To prove the theorem, suppose $\begin{pmatrix} 0 & 2 \\ 2 & 0 \end{pmatrix}$ embeds primitively in $\Pic{X}$.  Denote by $D, \Delta$ the classes spanning the $2 \times 2$ matrix. If the quasi-polarization $L$ on the surface $X$ is not pulled back from $S = \P^1 \times \P^1$, the Picard lattice of $X$ contains the $3 \times 3$ intersection matrix $M$ of the divisors $L, D, \Delta$, and $M$ is of the form 
\[ M = \begin{pmatrix}
    2g - 2 & a & b \\ a &0 & 2 \\ b & 2 & 0
\end{pmatrix}\]
where $a,b \ge 0$ and $\det M > 0$ by the Hodge index theorem.  These loci have codimension 2 and the requirement $\det M > 0$ holds for all $a,b \gg 0$ so it is straightforward to see there exist countably infinitely many components.  However, a priori, these codimension 2 components may be contained in divisorial components parametrizing K3s with $\irr(X) = 2$.

We first claim that $M$ does not represent $2$ if and only if $a,b$ are even and $g$ is odd.  Writing $V = xL + yD + z\Delta$, the equation $V^2 = 2$ halves to
\[(g - 1)x^2 + x(ay + bz) + 2yz = 1.\]
If $x$ is even, then the left side of the equation is even, so we cannot have it equal to $1$.  Hence, $x$ is odd, but then reducing modulo $2$ implies $ax + by \equiv g \pmod 2$, which is impossible if $g$ is odd and $a,b$ are even. Conversely, if $a$ or $b$ is odd, or $g$ is even, there is an explicit solution: 
\begin{itemize}
    \item if $a$ is odd, then $x = 1, y = 2 - g - \frac{b(1 - a)}{2}, z = \frac{1-a}{2}$ is a solution; 
    \item if $b$ is odd, then $x = 1, y = \frac{1-b}{2}, z = 2 - g - \frac{a(1 - b)}{2}$ is a solution; and 
    \item if $a,b$ are both even and $g$ is even, then $x = 1, y = 1 - \frac{b}{2}, z = \frac{ab/2 - g + 2 - a}{2}$ is a solution.
\end{itemize}

In any of the cases that $M$ does represent $2$, the intersection matrix of the quasi-polarization and the class of square $2$ is a Noether--Lefschetz divisor containing K3s with Picard lattice $M$, and by Lemma \ref{lem:class-of-square-2}, all K3s parametrized by this divisor have $\irr(X) = 2$.  If $M$ does not represent $2$, then $M$ cannot be contained in the Noether--Lefschetz divisors from Lemma \ref{lem:class-of-square-2} or $A = \begin{pmatrix} 0 & 2 \\ 2 & -2\end{pmatrix}$ or  $B= \begin{pmatrix} 0 & 1 \\ 1 & -2\end{pmatrix}$, as all of these represent 2. The only other Noether--Lefschetz divisor representing K3 surfaces with $\irr(X) = 2$ is that corresponding to the matrix $P = \begin{pmatrix} 0 & 2 \\ 2 & 0\end{pmatrix}$.  We claim that this matrix embeds in $M$ containing the polarization $L$ if and only if $a$ or $b$ is equal to $2$.  Given this claim, if $a,b$ are not equal to 2, then $M$ is not contained in any Noether--Lefschetz divisor representing K3s with $\irr(X) = 2$, and in this case the quasi-polarization $L$ is not pulled back from any map to a Hirzebruch surface. 

It now suffices to prove the claim.  Suppose $g$ is odd and $a,b$ are even.  If $a = 2$ (the case $b = 2$ is similar), then the class $D$ pairs with $L$ as $L \cdot D = 2$, and the lattice spanned by $L - \frac{g-1}{2}D$ and $D$ contains $L$ and has intersection matrix $P$.  For the reverse direction, write $a = 2a'$, $b = 2b'$, and $\det M = 4d$ where $d = ab - (2g-2)$.  Write $d = 2d'$ with $d' = 2a'b' - (g-1)$.  Observe that the class $u = L - b'D - a'\Delta$ pairs trivially with $D$ and $\Delta$ and has $u^2 = u \cdot L = -2d'$.  Consider the basis $\langle D, \Delta, u \rangle$ for the lattice $M$. 

Suppose there exist classes $E, F$ with intersection matrix $P$ embedding in $M$ such that $L \in \langle E, F \rangle$.  This must be primitive as $M$ does not represent 2.  The lattice $M$ must decompose as $\langle E, F \rangle \oplus \Zb v$ for a vector $v$ orthogonal to $E$ and $F$ with $v^2 = -2d'$ and $v \cdot w \in 2d'\Zb$ for any $w \in M$.  As $L \in \langle E, F \rangle$, we have $v \cdot L = 0$.  Write $v = x D + y \Delta + z u$ in the previous basis for $M$.  Then, 
\[ v \cdot D = 2y \in 2d'\Zb \implies y = d' y',\]
\[ v \cdot \Delta = 2x \in 2d'\Zb \implies x = d' x',\]
\[ v \cdot L = 2a' x + 2b' y - 2d' z =0 \implies z = a'x' + b' y', \] and
\[ v^2 = 4xy - dz^2 = -d \]
which, by substituting the third equation into the last and simplifying, implies 
\[ (a'x'+b'y')^2 - 2d'x'y' = 1.\]
If $x' y' < 0$ then the left hand side is at least $2d'|x'y'| \ge 2$, a contradiction.  If $x'y' > 0$, the left side is  
\[ a'^2 x'^2 + b'^2 y'^2 + (2a'b' - 2d')x'y' = (a'x' - b'y')^2 + (4a'b' - 2d')x'y'\]
and recall $d' = 2a'b' - (g-1)$, so this is at least $2(g-1)x'y' \ge 2$, a contradiction.  So, we must have $x'y' = 0$.  Without loss of generality, assume $y' = 0$, so the equation forces $a'^2 x'^2 = 1$, which implies $a' = 1$ and hence $a = 2$.  Therefore, we conclude $P = \begin{pmatrix} 0 & 2 \\ 2 & 0\end{pmatrix}$ embeds in $M$ containing the polarization $L$ if and only if $a$ or $b$ is equal to $2$.  In particular, for $a,b > 2$, $M$ cannot be contained in a Noether--Lefschetz divisor parametrizing K3s with $\irr(X) =2 $.

In summary, if $g$ is even, the locus of K3 surfaces with $\irr(X) = 2$ is a countably infinite union of Noether--Lefschetz divisors.  If $g$ is odd, the locus of K3 surfaces with $\irr(X) = 2$ is a countably infinite union of Noether--Lefschetz divisors and countably infinite codimension 2 loci corresponding to K3 surfaces for which the matrix \[ M = \begin{pmatrix}
    2g - 2 & a & b \\ a & 0 & 2 \\ b & 2 & 0
\end{pmatrix}\]
embeds into $\Pic{X}$, for some even integers $a,b \ge 4$.  
\end{proof}

By \Cref{ex:irrPicZK3}, we conclude the following.

\begin{corollary}\label{cor:non-constructible-K3s}
    For $g \geq 3$, the locus $\Hc_{g} \subset \Fc_{g}$ of K3 surfaces with $\irr = 2$ is not constructible. 
\end{corollary}

\section{Constant irrationality loci are not constructible in families}\label{section:IrratLociNotConstructible}

The goal of this section is to complete the proof of \Cref{cor:G2g} and provide another instance of Theorem \ref{thm:main-result} in certain families of general type surfaces. Recall the notation $\gonlocus{e}{g}$ and $\Hc_{g}$ from the introduction.

\begin{proof}[Proof of \Cref{cor:G2g}]
    Suppose $(X,L) \in \Gc_{2;g}$ for some $g \ge 3$.  By definition, $X$ contains a smooth curve of gonality $2$.  If this curve has genus 1, then $X$ admits an elliptic fibration and is contained in a Noether--Lefschetz divisor corresponding to $\Lambda = \begin{pmatrix} 2g-2 & a \\ a & 0 \end{pmatrix}$ for some $a > 0$.  Furthermore, every K3 surface with this Picard lattice contains a smooth elliptic curve. Now suppose this curve has genus at least 2, which by Theorem \ref{thm:irr2hypcurve} implies that $\irr(X) = 2$.  By Theorem \ref{thm:Irr2K3s}, $\Pic{X}$ either admits an ample class $A$ of square $2$ or one of the lattices in Theorem \ref{thm:Irr2K3s}(b) admits an embedding into $\Pic{X}$.  In the second case, $X$ contains an elliptic curve, so by the argument above it is contained in a Noether--Lefschetz divisor. In the first case, $X$ is contained in the Noether--Lefschetz divisor corresponding to $\Lambda = \begin{pmatrix} 2g-2 & a \\ a & 2 \end{pmatrix}$ for some $a > 0$, which is the intersection matrix of the quasi-polarization $L$ and the ample class $A$.  Furthermore, every quasi-polarized K3 surface $(X',L')$ parametrized by this Noether--Lefschetz divisor contains a class $A'$ of square $2$, which by Theorem \ref{thm:Irr2K3s-version2}, have $\irr(X') = 2$. This proves that $\gonlocus{2}{g}$ is a union of Noether--Lefschetz divisors in $\Fc_g$ for $g \ge 3$, and by Theorem \ref{thm:codim2-hyperelliptic}, this union is countably infinite. 
\end{proof}

In order to produce the families of surfaces of general type, we will first prove a lemma about double covers.

\begin{lemma}\label{lem:doubleCoverQuartics}
Let $Y$ be a smooth variety of dimension $n$ with a very ample line bundle $H$. Fix a smooth divisor $D \in |2aH|$ and let $f \colon X \rightarrow Y$ be the double cover branched over $D$. If $K_{Y} + aH$ separates $2 \cdot \irr(Y)$ points on $Y$, then $\irr(X) = 2 \cdot \irr(Y)$, and moreover every map $X \dashrightarrow \P^{n}$ of degree equal to $2 \cdot \irr(Y)$ factors through $f$.
\end{lemma}

\begin{proof}

Composing with $f$ implies that $\irr(X) \leq 2 \cdot \irr(Y)$. Now let $g \colon X \dashrightarrow \P^{n}$ be any generically finite map of degree $\leq 2 \cdot \irr(Y)$, and let $Z \colonequals g^{-1}(p)$ be the preimage of a general point $p \in \P^{n}$. Our goal will be to show that $g$ factors through $f$.

\textbf{Step 1:} We claim that the map $Z \rightarrow f(Z)$ is of constant degree. Suppose for contradiction that this is not the case. Let $U \subset \P^{n}$ be the open subset over which the map $g$ is \'{e}tale, and let $d$ be the largest integer for which a subscheme of $Z$ of length $d$ maps to a single point in $\P^{n}$ under $f$. Since $\# Z > \# f(Z)$ by assumption, it follows that $2 \leq d \leq \min \{ \deg g-1, \deg f \} = 2$ and hence $d = 2$. Consider the subset
\[ V =\{ [x_{1} + x_{2}] \in \Sym^{2}(X) \setminus \Delta \mid g(x_{1}) = g(x_{2}) \}. \]
Here $\Delta$ denotes the diagonal. Taking the preimage of $V$ under $X^{2} \rightarrow \Sym^{2}(X)$ and then projecting onto the first factor via $X^{2} \rightarrow X$, we get an algebraic subset $W$ of $g^{-1}(U)$ of dimension $n$ whose complement $g^{-1}(U) \setminus W$ is also algebraic of dimension $n$ (since both $W$ and $g^{-1}(U) \setminus W$ dominate $U$ under the map $g$). This contradicts the fact that $g^{-1}(U)$ is irreducible.

\textbf{Step 2:} Next we use trace maps to show that $f$ is not injective on $Z$. Assume for contradiction that this is the case, i.e. $\# Z = \# f(Z)$. Since $\# f(Z) = \# Z = \deg g \leq 2 \cdot \irr(Y)$, by assumption we can find a section $\eta \in H^{0}(Y, K_{Y} + aH)$ that vanishes at all but one of the points of $f(Z)$, and is nonzero at the last point. On the other hand, the Riemann--Hurwitz formula says that $K_{X} \cong f^{\ast} (K_{Y}+aH)$ so $f^{\ast}\eta \in H^{0}(X, K_{X})$ is a section that vanishes at all but one of the points of $Z$. Tracing this forward (see \cite{Mumford98}) via $g$ gives a nonzero section $\text{Tr}_{g}(f^{\ast}\eta) \in H^{0}(\P^{n}, K_{\P^{n}})$, which is a contradiction. Thus, we conclude that $f$ is not injective on $Z$.

\textbf{Step 3:} From the previous two steps, we conclude that $Z \rightarrow f(Z)$ must have constant degree equal to 2. So $\deg g = 2k$ for some positive integer $k$, and $Z$ is of the form
\[ Z = f^{-1}(\{ p_{1} \}) \sqcup \cdots \sqcup f^{-1}(\{ p_{k} \}) \]
for some points $p_{1}, \ldots, p_{k} \in Y$. Thus, by \cite[Lemma 1.3]{CM23} the map $g$ factors birationally into $h \circ f$ for some map $h \colon Y \dashrightarrow \P^{n}$ of degree $k$. This establishes the last claim of the lemma. But by definition $k \geq \irr(Y)$ and so $\deg g = 2k \geq 2 \, \irr(Y)$. Combining this with the first sentence in the proof yields the equality $\irr(X) = 2 \cdot \irr(Y)$.
\end{proof}

\begin{example}\label{ex:generaltype-non-constructible}
Let $\Sc \rightarrow T$ be the universal family of smooth quartic K3 surfaces. Let $A$ be the relative hyperplane class on $\Sc$, which is globally generated since $T$ is affine. Let $\Xc \rightarrow \Sc$ denote the double cover branched over a smooth divisor $D \in |2kA|$, for some integer $k \geq 5$. This gives a family $\pi \colon \Xc \rightarrow T$ whose general member is a smooth surface of general type. By generic smoothness, there exists an open subset $T^{\circ} \subset T$ such that the restricted family $\pi^{\circ} \colon \Xc^{\circ} \colonequals \pi^{-1}(T^{\circ}) \rightarrow T^{\circ}$ is a smooth family of general type surfaces. By \Cref{lem:doubleCoverQuartics}, the very general member of $\pi^{\circ}$ has degree of irrationality equal to 6, but by \Cref{thm:codim2-hyperelliptic} the locus in $T^{\circ}$ where the fibers have degree of irrationality equal to 4 is a countably infinite union of divisors.
\end{example}

Finally, in the families of surfaces constructed above, we claim that the locus of points in the base where the degree of irrationality drops is dense in the Euclidean topology. As noted by J. C. Ottem, this follows from the fact that K3 surfaces with Picard rank at least 12 admit an elliptic fibration with a section \cite[Corollary 14.3.8]{Huybrechts16}, and hence have $\irr(X) = 2$. In particular, it is enough to show that the Kummer locus is dense in $\Fc_{g}$, which follows from \cite[Remark 14.3.25]{Huybrechts16}.

Combining Corollary \ref{cor:non-constructible-K3s} and Example \ref{ex:generaltype-non-constructible} completes the proof of Theorem \ref{thm:main-result}.

\begin{remark}
In higher dimensions, one can construct families exhibiting similar behavior by combining \Cref{cor:non-constructible-K3s} with the following observation. For a surface $X$ with $\kappa(X) \geq 0$, one can show that $\irr(X \times \P^{m}) = 2$ if and only if $\irr(X) = 2$. Indeed, since $X$ is not covered in rational curves, any involution $\iota$ on $X \times \P^{m}$ must descend to an involution $\iota_X$ on $X$. The quotient $X / \iota_{X}$ is then unirational, and hence rational.
\end{remark}

\begin{remark}
    Using the main theorem of \cite{HJY26}, one can prove a nonconstructibility statement for the moduli spaces $\Ac_{1,d}$ of $(1,d)$ polarized abelian surfaces, when $d \not= 1, 2, 3, 6$. Namely, for these values of $d$ the locus of abelian surfaces in $\Ac_{1,d}$ that have $\irr = 3$ is not constructible.
\end{remark}

\section{Examples involving higher gonality curves on K3 surfaces}\label{section:irr3-and-trigonal}

In this section, we will use a result of Green--Lazarsfeld \cite{GL87} to construct an example of a K3 surface with degree of irrationality 3 that does not contain any smooth trigonal curves.

\begin{prop} \label{prop:picard_rank_1_clifford}
Suppose that $X$ is a K3 surface with $\Pic{X} = \Zb \cdot L$ for some ample $L$. Then any smooth curve $C \in |L|$ has Clifford index equal to its generic value $\Cliff(C) = \floor{(g-1)/2} = \floor{L^2/4}$.
Furthermore, any smooth curve $C \in |m L|$ with $m > 1$ has Clifford index $\Cliff(C) = L^2 (m-1) - 2$,
which is strictly less than the generic value for curves of its genus.
\end{prop}

\begin{proof}
By the main theorem of Green--Lazarsfeld \cite{GL87} (see \Cref{thm:GL} for the statement), for any fixed $m \ge 1$ every smooth curve in the linear series $C \in |mL|$ has the same Clifford index. Moreover, if $\Cliff(C)$ is strictly less than the generic value $\lfloor \frac{g-1}{2} \rfloor$, then there is a line bundle on $X$ whose restriction computes the Clifford index of $C$.

Let $C \in |mL|$ be a smooth curve. By Kodaira vanishing, we get an exact sequence $0 \to H^0(\Oc_X((\ell - m) L)) \to H^0(\Oc_X(\ell L)) \to H^0(\Oc_C(\ell L)) \to 0$. One can then use this to compute the cohomology of the restriction to $C$ of the line bundle $\ell L$ on $X$, for some $\ell \in \Zb$. It follows that the cases where $l \ge m$ or $l \le 0$ do not satisfy both of $h^0(\Oc_C(\ell L)) \ge 2$ or $h^1(\Oc_C(\ell L)) \ge 2$, and so by definition they cannot compute the Clifford index of $C$.

In the remaining cases $0 < \ell < m$, one computes that 
\[ \Cliff(\Oc_C(\ell L)) = m\ell L^2 - 2 (h^0(\Oc_C(\ell L)) - 1) = L^2 ( m \ell - \ell^2) - 2. \]
For $m = 1$ it follows that $\Cliff(C)$ cannot be computed by the restriction of a line bundle on $X$. Thus, by \cite{GL87} it must achieve the generic value for $m = 1$. For $m > 1$, the smallest value of $\Cliff(\Oc_C(\ell L))$ occurs at either $\ell = 1$ or $\ell = m-1$, which give the same value. In either case, we see that
\[ \Cliff(\Oc_C(L)) = \Cliff(\Oc_C((m-1)L)) = L^2 (m-1) - 2 \]
and it is straightforward to check that $L^2 (m-1) - 2 <  \floor{m^2 L^2 / 4} =\floor{(g-1)/2}$ when $m \geq 2$. Since this implies that $\Cliff(C)$ is strictly less than the generic value, by \cite{GL87} there is some line bundle on $X$ whose restriction to $C$ computes $\Cliff(C)$. From the computation of $\Cliff(\Oc_{C}(\ell L))$ above, it follows that this line bundle must be either $\Oc_{C}(L)$ or $\Oc_{C}((m-1)L)$.
\end{proof}

We will use the proposition above to prove:

\begin{theorem}\label{thm:irr3-but-no-trigonal}
Let $(X, L)$ be a polarized K3 surface of genus $6$ with $\Pic{X} = \Zb \cdot L$. Then $\irr(X) = 3$ but $X$ does not contain a smooth trigonal curve.
\end{theorem}

\begin{proof}
Let $C \subset X$ be a curve, which necessarily lies in $|mL|$ for some $m > 0$. Recall  that $\Cliff(C) = 1$ if and only if $C$ is either a trigonal curve or a smooth plane quintic. Therefore, it suffices to show that $C$ does not have Clifford index $1$. By Proposition \ref{prop:picard_rank_1_clifford}, if $m = 1$ then $\Cliff(C) = 2 \ne 1$, and if $m > 1$ then $\Cliff(C) \ge 10(m-1) - 2 \ge 8 \ne 1$. To complete the proof, note that \cite[Theorem A]{Moretti23} implies that $\irr(X) \le 3$ but \Cref{ex:irrPicZK3} gives $\irr(X) > 2$.
\end{proof}

\section{Polarized gonality loci on K3 surfaces}\label{section:polarized-gonality-loci}

The main goal of this section is to prove that polarized gonality loci are divisors in $\Fc_{g}$ (see \Cref{thm:polLociDivisors}). Recall from the introduction the definition of the unpolarized $e$-gonal loci $\gonlocus{e}{g}$. One can similarly define $\gonlocus{\leq e}{g}$ where $C$ has gonality at most $e$.

We may refine these definitions by allowing $C$ to lie in some multiple of the polarization.

\begin{definition}
Let $\Fc_g$ be the moduli space of quasi-polarized K3 surfaces of genus $g$. Let $e > 0$ be an integer. Then we define the \textit{polarized $e$-gonal locus} of $\Fc_g$ as:
\begin{enumerate}
    \item $\pgonlocus{e}{g}^{m} \colonequals \{ (X, L) \in \Fc_{g} \colon \exists \text{ smooth curve } C \subset X \text{ with } C \in |mL| \text{ and } \gon(C) = e \}$
    \item $\pgonlocus{e}{g} \colonequals \bigcup_{m > 0} \pgonlocus{e}{g}^{m}$.
\end{enumerate}
Similarly, the \textit{polarized $\le e$-gonal locus} is
\begin{enumerate}
    \item $\pgonlocus{\le e}{g}^{m} \colonequals \{ (X, L) \in \Fc_{g} \colon \exists \text{ smooth curve } C \subset X \text{ with } C \in |mL| \text{ and } \gon(C) \le e \}$
    \item $\pgonlocus{\le e}{g} \colonequals \bigcup_{m > 0} \pgonlocus{\le e}{g}^{m}$.
\end{enumerate}
\end{definition}

To avoid confusion, we will refer to the loci $\gonlocus{e}{g}$ and $\gonlocus{\le e}{g}$ as the \textit{unpolarized gonality loci} in what follows. Note that $\Gc_{1;g}$ is simply the locus of K3 surfaces containing a smooth rational curve and that $\pgonlocus{1}{g}$ is always empty.

\begin{proposition}
    The loci $\gonlocus{\leq e}{g}$ and $\pgonlocus{\le e}{g}$ in $\Fc_{g}$ are closed under specialization.
\end{proposition}

\begin{proof}
    By \cite[Chapter 17, Proposition 2.10]{Huybrechts16}, the specialization map on $\Pic$ is injective in a smooth family of quasi-polarized K3 surfaces and the intersection pairing is preserved. Thus, it suffices to assume that $e \geq 2$. Let $C \in |M|$ be a smooth curve of gonality $\leq e$ on the very general fiber $X$ of a 1-parameter family of K3 surfaces (equivalently the geometric generic fiber), and suppose that $M$ specializes to some line bundle $M'$ on $X'$.
    
    In the polarized case $\pgonlocus{\le e}{g}^{m}$, the line bundle $M$ is a multiple $mL$ of the quasi-polarization on $X$, and so $M'$ is automatically big and nef.  Since we are working with a family in $\Fc_{g}$, the line bundle $L$ specializes to some $L'$ which is also primitive. If $M'$ is not basepoint-free, then $M'$ is of the form $M' = mE + R$ as in \Cref{prop:base-locus-of-big-and-nef}(2) and $eM'$ is basepoint-free when $e \geq 2$. A general section of $eM'$ is a degree $e$ multisection of the elliptic fibration $|E|$. Thus, this multisection has gonality at most $e$, and the special fiber $X'$ lies in $\pgonlocus{\le e}{g}$.
    
    In the unpolarized case $\gonlocus{\leq e}{g}$, if $M$ is not a multiple of the quasi-polarization (which happens when we consider a family in $\gonlocus{\leq e}{g}$), then it is possible for the specialization $M'$ of $M$ to fail to be nef. But in this case, the denefing curve is a smooth rational curve and so the special fiber lies in $G_{1;g}$ and we are done. If $M'$ fails to be basepoint-free, then by \Cref{prop:base-locus-of-big-and-nef}(2) there is a smooth rational curve on $X'$ and so we are again done.
    
    Thus, from now on we may assume that $M$ specializes to a line bundle $M'$ on $X'$ that is basepoint-free. If $|M|$ is not equal to the Donagi--Morrison example \cite[(2.2)]{DonagiMorrison}, then by \cite[Theorem 2]{Knutsen09} the gonality is constant for any smooth member of the linear series $|M|$. The entire linear series $|M|$ specializes to $|M'|$ since we have the vanishing $h^{1}(\Oc_{X'}(M')) = 0$. Replacing $C$ with a general member of $|M|$ (which has the same gonality) and specializing, it follows that any smooth curve $C' \in |M'|$ satisfies $\gon(C') \leq \gon(C)$.
    
    Now suppose $|M|$ is the Donagi--Morrison example. Recall from \Cref{ex:DM} that this means that $\pi \colon X \rightarrow \P^{2}$ is a double cover of $\P^{2}$ branched along a sextic and $M = \pi^{\ast}\Oc_{\P^{2}}(3)$ is a genus 10 polarization. The smooth curves in the codimension one linear subspace $|\pi^{\ast}H^{0}(\Oc_{\P^{2}}(3))| \subset |M|$ are bielliptic and have gonality 4, whereas the general smooth curve in $|M|$ has gonality 5. Since the Clifford index has to be constant by \cite{GL87}, these are the only possibilities for the gonality. If $C \in |M|$ has gonality 5, then the same argument as in the previous paragraph applies. If $C \in |M|$ has gonality 4, then in the polarized case the linear series $|L'|$ defines a map from $\pi' \colon X' \rightarrow \P^{2}$ of degree 2 (using \Cref{prop:base-locus-of-big-and-nef} and the fact that $(L')^{2} = 2$). A general plane cubic will meet the branch divisor of $\pi'$ transversely, and so its preimage is a smooth curve of gonality four and we are done. In the unpolarized case, we observe that $|L'|$ is big and nef of degree 2, so either there exists a smooth rational curve on $X'$ or $|L'|$ is basepoint-free and contains a smooth hyperelliptic curve. In either case, we are done.
\end{proof}

Next, we give a basic structural result relating the unpolarized gonality locus $\gonlocus{e}{g}$ to the polarized gonality loci. Suppose $(X,L) \in \Fc_{g}$ is a quasi-polarized K3 surface in $\gonlocus{e}{g}$. From the definition of $\gonlocus{e}{g}$, there exists a smooth integral curve $C \subset X$ with $\gon(C) = e > 2$ and $C^2 \ge 2$, so $\Oc(C)$ is basepoint-free \cite[Chapter 2, Proposition 3.5]{Huybrechts16}. Hence, writing $\Oc(C) \sim (L')^{\otimes m}$ for $m > 0$ maximal we see that $(X, L')$ is quasi-polarized.  As $L'^2 = 2 h - 2$, we have $(X, L') \in \mathcal{F}_{h}$.    If $L' \ne L$, by definition, $X$ carries two quasi-polarizations $L$ and $L'$. Let $\Lambda = \langle L, L' \rangle$ and consider the Noether--Lefschetz divisor $\NL_{\Lambda; h} \subset \Fc_h$ of quasi-polarized K3 surfaces $(Y,L'_Y)$ of genus $h$ such that $\Lambda \subset \Pic{Y}$.  Write $\Pol^{\Lambda,L}_{h} \subset \NL_{\Lambda;h} \subset \mathcal{F}_h$ for the locus such that the class $L_Y$ on $Y$ remains big and nef.  If $L' = L$, then define $\Pol^{\Lambda,L}_{g} = \Fc_g$.

For $L' \ne L$, observe that $\Pol^{\Lambda,L}_{h}$ is an open subset of $\NL_{\Lambda; h}$.  Indeed, assuming $L$ is big and nef, we have $H^i(X,L^{\ot m}) = 0$ for $i > 0$ and $m > 0$ by Kawamata-Viehweg vanishing. By upper semi-continuity, the same holds for any $(Y, L'_Y)$ in an open neighborhood of the point $(X, L') \in \NL_{\Lambda; h}$.  On this open set, the number of global sections of $H^0(Y, L_Y^{\ot m})$ are constant so we conclude there is an open set on which $L$ remains big.  As nef implies semiample for K3 surfaces, the line bundle $L$ on $X$ is semiample, so for $m \gg 0$ the map $H^0(X,L^{\ot m}) \otimes \Oc_X \to L^{\ot m}$ is surjective.   Writing $\pi: (\Xc, \Lc) \to \NL_{\Lambda}$ for the universal family over $\NL_{\Lambda}$ and using cohomology and base change, this means $\pi^*\pi_*\Lc^{\ot m} \to \Lc^{\ot m}$ is surjective at the point $(X,L)$ and hence is surjective in a neighborhood.  Therefore, there is an open neighborhood of $(X,L)$ such that $L_Y$ is semiample and hence nef for all $(Y,L_Y)$ in the neighborhood.  We conclude there is an open neighborhood $U$ of $(X,L') \in\Pol^{\Lambda,L}_{h}$ such that for all $(Y,L'_Y) \in U$, $L_Y$ remains big and nef.  Hence, $\Pol^{\Lambda,L}_{h}$ is open in $\NL_{\Lambda; h}$.

The argument in the previous paragraph now proves the formula
\[ \Gc_{e;g} = \bigcup_{\substack{h, \Lambda \\ L \in \Lambda : L^2 = 2h - 2}} \im{(\pgonlocus{e}{h} \cap \Pol^{\Lambda,L}_{h} \to \Fc_g)} \]
where the natural map $\Pol^{\Lambda, L}_{h} \to \Fc_g$ is defined by using the quasi-polarization $L$. Note that the sets $\Pol^{\Lambda,L}_{h}$ are locally closed in $\Fc_h$. With this formula, we prove the following:

\begin{prop}
$\pgonlocus{e}{g}^{m}$ is locally closed. Hence, $\pgonlocus{e}{g}$ and $\gonlocus{e}{g}$ are countable unions of locally closed subsets. 
\end{prop}

\begin{proof}
$\pgonlocus{e}{g}^{m}$ is locally closed by Chevalley's theorem and constructibility of the gonality $e$ locus in a family of curves. Then the identifications
\[ \pgonlocus{e}{g} \colonequals \bigcup_{m > 0}\pgonlocus{e}{g}^{m} \qquad \text{and} \qquad \gonlocus{e}{g} = \bigcup_{\substack{h, \Lambda \\ L \in \Lambda : L^2 = 2h - 2}} \im{(\pgonlocus{e}{h} \cap \Pol^{\Lambda, L}_{h} \to \Fc_g)} \]
together with the local closedness of $\Pol^{\Lambda, L}_{h}$ in $\Fc_h$ prove the claim. 
\end{proof}

In the remainder of this section, we would like to relate the gonality loci to Noether--Lefschetz divisors. The main idea is to pass to the Clifford index, where the Green--Lazarsfeld type results apply, and use the fact that the Clifford dimension is a quantitative measure of the discrepancy between the Clifford index and the gonality. For this, we will need two preparatory lemmas.

\begin{lemma}\label{lem:prop-case}
    Let $(X,L)$ be a quasi-polarized K3 surface and let $C \in |mL|$ be a non-hyperelliptic curve with non-generic gonality for its genus, for some $m > 0$.  Let $D, \Delta \in \Pic{X}$ be the divisors in Theorem \ref{thm:DonagiMorrison}.  If $D, \Delta$ are proportional, then: 
    \begin{enumerate}
        \item Up to interchanging $D, \Delta$, we have $D \sim L$ and $\Delta \sim n L$ for some $n > 0$, and $\Cliff(C) = nL^2 - 2$.
        \item If $C$ is not exceptional, then $\gon(C) = nL^2$.
        \item If $C$ is exceptional, then $(X,L)$ and $C$ are as in Example \ref{ex:DM}.
    \end{enumerate}
\end{lemma}

\begin{proof}
Note that Theorem \ref{thm:DonagiMorrison} applies as $C$ is not hyperelliptic and has non-generic gonality for its genus.  By assumption, $D, \Delta$ are proportional.  Since $C \in |D + \Delta|$ and $C \in |mL|$, it follows from the primitivity of $L$ that $D \sim n_1L$ and $\Delta \sim n_2 L$ for some $n_1, n_2 > 0$.  Because $L$ is a quasi-polarization, we also know that $D, \Delta$ are big and nef.  From the exact sequence
\[ 0 \to \Oc(-\Delta) \to \Oc(D) \to \Oc_C(D|_C) \to 0, \]
it is straightforward to compute $\Cliff(D|_C) = n_1n_2L^2 - 2$.   

We claim that $n_1 = 1$. From Remark \ref{rmk:DM-g1d}, if $C$ is not exceptional then $\Cliff(C) = \Cliff(D|_C)$, and if $C$ is exceptional, we have $\Cliff(D|_C) -1 \le \Cliff(C)$.  From the exact sequence 
\[ 0 \to \Oc(L-C) = \Oc(- (n_1 + n_2 - 1)L) \to \Oc(L) \to \Oc_C(L|_C) \to 0 \] observe that $L$ contributes to $\Cliff(C)$ because $h^0(C, \Oc(L|_C)) = h^0(X, \Oc(L)) = 2 + \frac{L^2}{2}$ and $h^1(C, \Oc(L|_C)) = h^0(X, \Oc((n_1 + n_2 - 1)L)) \ge 2$. Because $\Cliff(L|_C) = (n_1 + n_2 - 1)L^2 - 2$, and in either case we must have $\Cliff(D|_C) -1 \le \Cliff(C) \le \Cliff(L|_C)$, this implies 
\[ n_1n_2L^2 - 3 \le (n_1 + n_2 - 1)L^2 - 2\] or equivalently $(n_1 -1)(n_2 - 1)L^2 \le 1$. As $L^2> 0$ is even, this is possible if and only if either $n_1 = 1$ or $n_2 = 1$.  Interchanging the roles of $D$ and $\Delta$ if necessary, we conclude $L \sim D$, $\Delta \in |n D|$ for some $n > 0$, and $C \in |(n+1) D|$ and $D \cdot \Delta = nL^2$.

If $C$ is not exceptional, we then conclude $\Cliff(C) = \Cliff(D|_C) = n L^2 - 2$ and hence $\gon(C) = nL^2$.  If $C$ is exceptional, either $\Cliff(C) = nL^2 - 2$ or $\Cliff(C) = nL^2 - 3$. We claim that $(X,C)$ cannot be a generalized ELMS example and this will rule out the second case.  From above, we have $C \in |(n+1) D|$, and if $(X,C)$ were a generalized ELMS example, then by definition there exists some $D', \Gamma$ on $X$ such that $(D')^2 \ge 2$, $\Gamma^2 = -2$, and $D' \cdot \Gamma = 1$ such that $C \in |2D' + \Gamma|$.  This implies $(n+1)D \sim 2D' + \Gamma$ or equivalently $(n+1)D - 2D' \sim \Gamma$.  Intersecting both sides with $\Gamma$ and using that $D' \cdot \Gamma = 1$ implies $D \cdot \Gamma = 0$.  This is sufficient to compute $(D')^2$: 
\[(n+1)L^2 = D \cdot C = D \cdot (2D' + \Gamma)\]
implies $D \cdot D' = \frac{n+1}{2} L^2$, and then 
\[\tfrac{(n+1)^2}{2} L^2 = (n+1) D \cdot D' = C \cdot D' = (2D' + \Gamma) \cdot D' = 2(D')^2 + 1 \]
implies $(D')^2 = \frac{(n+1)^2 L^2 -2}{4}$.  By definition of the generalized ELMS examples, $\Cliff(C) = (D')^2 -1 =  \frac{(n+1)^2 L^2 -6}{4}$, and from above, we must also have $\Cliff(C) = nL^2 - 2$ or $nL^2 - 3$. In either case, since $L^{2}, n > 0$ we get a contradiction. 

Thus, $(X,C)$ cannot be a generalized ELMS example. By \Cref{thm:knutsen} we conclude that $(X,C)$ is a Donagi--Morrison example, with $C \in |3L|$ of gonality 5 and $L^{2} = 2$, so $n = 2$. The curves in the Donagi--Morrison example have $\Cliff(C) = 2$ which is precisely $n L^2 - 2$ in this case. 
\end{proof}

The second lemma is about specializations of generalized ELMS examples.

\begin{lemma}\label{lem:ELMS-specialization}
Fix an integer $e \ge 3$ and let $X$ be a K3 surface with a sublattice $\Lambda = \langle D, \Delta \rangle^{\sat} \subset \Pic{X}$ with intersection pairing $\begin{pmatrix} e-1 & e \\ e & e-1 \end{pmatrix}$. Suppose that there exists a smooth curve $C \sim D + \Delta$ of genus $g = 2e$ with $\Cliff(C) = e - 2$. Then $(X,\Oc(C))$ is a generalized ELMS example. In particular, $\gon(C) = e + 1$ and $C$ is exceptional. 
\end{lemma} 

\begin{proof}
Note that $\Oc(C)$ is big and nef on $X$. To prove that $(X,\Oc(C))$ is a generalized ELMS example, we verify the definition. 

Our first goal will be to produce the desired smooth rational curve $\Gamma$. As a candidate, define $\Gamma := \Delta-D = C - 2D$ and observe that $\Gamma^2 = -2$ and $\Gamma \cdot C = 0$. Because $\Gamma^2 = -2$, by Riemann--Roch either $\Gamma$ or $-\Gamma$ is effective. If $\Gamma$ is effective, we have $D \cdot \Gamma = 1$. If $-\Gamma$ is effective, define $\Gamma' = - \Gamma$, $D' = D + \Gamma$ and $\Delta'= \Delta - \Gamma$. The intersection pairing of $D'$ and $\Delta'$ is then the same as that of $D$ and $\Delta$, and we still have $C \sim D' + \Delta'$ and $\Gamma' = \Delta' - D' = C - 2D'$. Thus if $\Gamma'$ is effective, we also have $D' \cdot \Gamma' = 1$. By the previous sentences, we may therefore assume that $\Gamma$ is effective. Since $C$ is big and semiample \cite[Chapter 2, Corollary 3.12]{Huybrechts16}, the intersection $\Gamma \cdot C = 0$ implies that $\Gamma$ is a disjoint union of ADE configurations of $(-2)$-curves \cite[Proposition 2]{Mayer}. Moreover, $\Gamma$ is connected because $\Gamma^2 = -2$. We will now show that $\Gamma$ consists of a single $(-2)$-curve.  If $R \subset \Gamma$ is an irreducible component of $\Gamma$, then
$R \cdot C = 0$ and $R^2 = -2$.  We claim that any $(-2)$-curve $R$ satisfying $R \cdot C = 0$ must either be equal to $\Gamma$ or be disjoint from it, from which it follows that $\Gamma$ is necessarily a single $(-2)$-curve. Indeed, as $\Gamma = C - 2D$ and $R \cdot C = 0$, the intersection number
\[ R \cdot \Gamma = R \cdot C - 2 R \cdot D = -2 R \cdot D \] 
is always even. If $R$ is a component of $\Gamma$, then $R \cdot \Gamma = -2 + \ell$ where $\ell$ is the number of neighbors of $R$ in the dual graph of $\Gamma$, which is an ADE configuration of curves and hence a tree.  Thus, $\ell$ must be even for each $R$. But on the other hand, each leaf of a tree has only one neighbor and therefore $\ell$ cannot be even for each component $R$ of $\Gamma$ unless $R = \Gamma$. If $R$ is not a component of $\Gamma$, then $(R + \Gamma)^2 = -4 + 2 (R \cdot \Gamma)$.  If $R$ meets $\Gamma$, then $(R+\Gamma)^2 \ge 0$ since $R \cdot \Gamma$ is even and positive.  However, $R + \Gamma$ is contracted by $|mC|$, and by negative-definiteness its square must be negative, so we conclude $R$ and $\Gamma$ are disjoint. This proves the desired claim. 
\par 
Next, we wish to show $D = \tfrac{1}{2}(C - \Gamma)$ can be chosen to be smooth.  First, we claim that $D$ is nef. Otherwise, by \cite[Corollary 2.1.5]{Huybrechts16} there must be some $(-2)$-curve $R'$ with $D \cdot R' < 0$. Consider the intersection matrix of the curves $C, \Gamma, R'$ 

\[
\begin{pmatrix}
4e-2 & 0 & C \cdot R'
\\
0 & -2 & \Gamma \cdot R'
\\
C \cdot R' & \Gamma \cdot R' & -2
\end{pmatrix},
\]
whose determinant
\[ (4 - (\Gamma \cdot R')^2) (4e-2) + 2 (C \cdot R')^2 > 0  \]
must be positive by the Hodge index theorem.  As $D \cdot R' < 0$, then $C \cdot R' < \Gamma \cdot R'$.  For ease of notation, write $\delta := -(D \cdot R')$ and $d := C \cdot R'$. Then $\Gamma \cdot R' = 2 \delta + d$, so the determinant becomes
\[ (4 - 4 \delta^2 - 4 \delta d - d^2) (4e-2) + 2 d^2 > 0, \]
and since $4e-2 \ge 2$, this implies
\[ 4 - 4 \delta^2 - 4 \delta d > 0 .\]
Note that $d \ge 0$ by nefness of $C$, but then this equation has no solutions with $\delta > 0$, which is a contradiction. Therefore, $D$ is nef, and $D^{2} = \frac{1}{4} (4e-4) = e-1 \geq 2$ so $D$ is big.

Now suppose that $D$ is not basepoint-free. By \Cref{prop:base-locus-of-big-and-nef}(2), it follows that $D \sim kE + R$ where $E$ is a smooth elliptic curve, $R$ is a smooth rational curve, $k \ge 2$ is an integer, and $E \cdot R = 1$. As $D \cdot \Gamma = 1$ and $E$ is nef, we see that $E \cdot \Gamma = \tfrac{1}{k} (D - R) \cdot \Gamma = \tfrac{1}{k}(1 - R \cdot \Gamma) \ge 0$ is a non-negative integer. Because $\Gamma$ is a smooth rational curve, we must have $R \cdot \Gamma \ge 0$ or $R = \Gamma$, so the possible cases are (1) $R \cdot \Gamma = 1$ and $k \ge 2$, or (2) $R = \Gamma$ and $k = 3$.

In the first case, we have $E \cdot \Gamma = 0$ and thus
\[ E \cdot C = 2 E \cdot D + E \cdot \Gamma = 2. \]
Since $E$ moves in a pencil, its restriction to $C$ shows that $C$ is hyperelliptic so $\Cliff(C) = 0$, contradicting our assumption that $\Cliff(C) = e - 2 > 0$. In the second case, $E \cdot \Gamma = 1$ so $D = 3E + \Gamma$ and therefore $D^2 = 6 E \cdot \Gamma - 2 = 4$. Therefore, we must have $e = 5$ but then $E \cdot C = 3$ implies that $\gon(C) \le 3$ so $\Cliff(C) \le 1$, contradicting the assumption that $\Cliff(C) = e-2 = 3$.  Therefore, we may assume $D$ is basepoint-free and smooth by Bertini.

Thus far, we have shown that $C \sim 2D + \Gamma$, where $D$ and $\Gamma$ are smooth curves satisfying $D^2 = e - 1\ge 2$, $\Gamma^2 = -2$, and $D \cdot \Gamma = 1$. To complete the proof, suppose that there exists $B$ with $0 \le B^2 < D^2$ and $0 < B \cdot C - B^2 \le D^2$.  We claim that this forces $\Cliff(C) < e - 2$ and so no such $B$ exists under our assumptions.  To prove the claim, we need to show that $B$ and $C - B$ are both effective with at least two sections.  Observe that $B^2 \ge 0$ and
\[ (C - B)^2 = C^2 - 2(B \cdot C - B^2) - B^2 \ge C^2 - 2 D^2 - B^2 > 2 D^2 - B^2 > 0,\] where the last inequality holds because 
\[ C^2 - 2D^2 = 2D^2 + 4D \cdot \Gamma + \Gamma^2 = 2D^2 + 2 .\]
Therefore, Riemann--Roch shows that, if they are effective, they have at least two sections.  But, $B \cdot C > B^2 \ge 0$ so $B$ is effective and
\[ (C - B) \cdot C = C^2 - (B \cdot C - B^2) - B^2 \ge C^2 - D^2 - B^2 >  D^2 - B^2 > 0,\] where the second to last inequality holds because 
\[ C^2 - D^2 = (2D + \Gamma)^2 - D^2 = 3D^2 + 4D \cdot \Gamma + \Gamma^2 = 3D^2 + 2 \]
so $C - B$ is also effective.  Therefore, $B|_C$ contributes to the Clifford index and 
\[ \Cliff(C) \le \deg{B|_C} - 2(h^0(B|_C) - 1) \le B \cdot C - 2 (\tfrac{B^2}{2} + 1) \le D^2 - 2 = e - 3. \]
By assumption, $\Cliff(C) = e - 2$, and hence no such $B$ exists. Thus, the pair $(X, \Oc(C))$ is a generalized ELMS example.
\end{proof}

Using these lemmas and the classification of K3 surfaces containing exceptional curves (from \S\ref{section:prelims}), we are now ready to prove Theorem \ref{intro-thm:pol-gonloci-divisors}.

\begin{theorem}\label{thm:polLociDivisors}
For fixed integers $e,m \geq 1$ and $g \geq 2$, the locus $\ol{\pgonlocus{e}{g}^{m}}$ is either empty, a finite union of Noether--Lefschetz divisors in $\Fc_g$, or is all of $\Fc_g$.
\end{theorem}

\begin{proof}
When $e = 1$, by definition $\ol{\pgonlocus{1}{g}^{m}}$ is empty. When $e=2$, the content of this theorem is contained in \cite{Reid76a} and Section \ref{section:hyperellipticK3s}, so we assume in what follows that $e > 2$.

Let $e > 2$ and assume $\ol{\pgonlocus{e}{g}^{m}}$ is non-empty. It suffices to show that either $\pgonlocus{e}{g}^{m}$ is dense, or $\ol{\pgonlocus{e}{g}^{m}}$ is contained in a finite union of Noether--Lefschetz divisors $B_1, \dots, B_k$ such that the generic point of each $B_i$ is contained in $\ol{\pgonlocus{e}{g}^{m}}$. Suppose that $(X, L) \in \pgonlocus{e}{g}^{m}$ and let $C \in |m L|$ witness $\gon(C) = e$. If $e = \floor{\frac{g(C) + 3}{2}}$ is maximal with respect to the genus of $C$, then we claim that the generic point of $\Fc_g$ is an element of $\pgonlocus{e}{g}^{m}$ and therefore $\ol{\pgonlocus{e}{g}^{m}} = \Fc_g$.  Indeed, for any deformation $(X',L') \in \Fc_g$ of $(X,L)$ and for $C' \in |mL'|$ a general member, $C'$ is a deformation of $C$ and therefore $g(C') = g(C)$. Because gonality can only decrease under specialization and $C$ had the maximal gonality for its genus, we conclude that $\gon(C') = e$ and therefore $(X',L') \in \pgonlocus{e}{g}^{m}$.

Now assume $2 < e < \floor{\frac{g(C) + 3}{2}}$. By Theorem \ref{thm:DonagiMorrison} and Corollary \ref{cor:DMmatrices}, $\Pic{X}$ contains effective classes $D$, $\Delta$ whose intersection pairing is given by one of the finitely many $(g(C), e)$ Donagi--Morrison matrices. 

\textbf{Step 1:} First assume $D, \Delta$ are proportional.  Then, by Lemma \ref{lem:prop-case}, we must have $D \sim L$, $\Delta \sim nL$, and $C \sim (n+1)L$.  Furthermore, if $C$ is exceptional, then $(X,L)$ must be as in the Donagi--Morrison example, which implies $(X, L)$ is a polarized K3 surface of genus $g = 2$,  $C \in |3 L|$ (so $m = 3$) and necessarily $e = 5$. In this case, the generic $(X', L') \in \Fc_2$ is a double cover of $\P^2$ branched over a smooth sextic and is a Donagi--Morrison example itself. Thus $|3 L'|$ also contains smooth curves of gonality $4$ and $5$. Hence, both gonality loci $\pgonlocus{4}{2}^{3}$ and $\pgonlocus{5}{2}^{3}$ are dense in $\Fc_2$.  

Now suppose $D, \Delta$ are proportional but $C$ is not exceptional. Because $D, \Delta, C$ are all proportional to the quasi-polarization, for very general $(X',L') \in \Fc_g$, there exist deformations $D', \Delta', C'$ on $X'$ satisfying the same proportionality.  If $C'$ is not exceptional, we have
\[D \cdot \Delta = \gon(C) \le \gon(C') \le D' \cdot \Delta' \] where the second inequality follows as $C$ is a specialization of $C'$ and the third follows from Lemma \ref{lem:DM-matrices-force-gon(C)}.  Because $D' \cdot \Delta' = D \cdot \Delta$, all inequalities must be equalities so we conclude the generic $(X',L') \in \Fc_g$ contains a curve of gonality $\gon(C)$, which means $\pgonlocus{e}{g}^{m}$ is dense.  By Lemma \ref{lem:prop-case}, if $C'$ is exceptional, we must have $(X',L')$ as in the Donagi--Morrison example with  $C'$ a smooth plane sextic with gonality $5$, $g = 2$, and $m = 3$.  From Lemma \ref{lem:prop-case}, we have  $\Cliff(C) = \Cliff(C') =2$ and $e = \gon(C) = 4$.  However, in this case, $X'$ also contains curves of gonality 4, so we see $\pgonlocus{4}{2}^{3}$ is dense in $\Fc_2$.

\textbf{Step 2:} Now suppose $D, \Delta$ are not proportional for any of the Donagi--Morrison matrices afforded by Theorem \ref{thm:DonagiMorrison} and Corollary \ref{cor:DMmatrices}, so $\Pic{X}$ contains a $2 \times 2$ sublattice spanned by $D, \Delta$. Define $B_1, \dots, B_k$ to be the set of Noether--Lefschetz divisors $\NL_{\Lambda;g}$ whose generic member is contained in $\pgonlocus{e}{g}^{m}$.  It suffices to show that $\pgonlocus{e}{g}^{m} \subset B_1 \cup \cdots \cup B_k$.

As $\gon(C) = e < \floor{\frac{g(C) + 3}{2}}$, from Theorem \ref{thm:DonagiMorrison} and Corollary \ref{cor:DMmatrices} there exists a Donagi--Morrison matrix $\Lambda = \left< D, \Delta \right>^{\sat}$ spanned by divisors with $\Cliff(D|_C) \le D \cdot \Delta - 2 \le e - 2 = \Cliff(Z)$, where $|Z|$ is the complete linear series corresponding to the $\g^1_e$ on $C$.  This chain of inequalities is proven in \cite[Lemma 4.6]{DonagiMorrison}.  We must show that the very general $(X', L') \in \NL_{\Lambda;g}$ is an element of $\pgonlocus{e}{g}^{m}$. We proceed in two cases based on the exceptionality of $C$.

\textbf{Step 2(a):} Suppose $C$ is not exceptional. Then $\Cliff(C) = \gon(C) - 2 = e - 2$ and therefore $\Cliff(C) = \Cliff(Z) = \Cliff(D|_C)$ and $D \cdot \Delta = e$.  By \cite[Proof of Lemma 4.6]{DonagiMorrison}, we must have $H^0(X,
\Oc(D)) \cong H^0(C,\Oc_C(D|_C))$ and $H^1(X,\Oc(D)) = 0$.  From the exact sequence 
\[ 0 \to \Oc(-\Delta) \to \Oc(D) \to \Oc_C(D|_C) \to 0\]
this also implies $H^1(X, \Oc(-\Delta)) = 0$, which by Serre Duality implies $H^1(X, \Oc(\Delta)) = 0$.

Let $(X',L')$ be very general in $\NL_{\Lambda;g}$ and let $C' \in |mL'|$ a smooth member.  We claim that $\Cliff(D|_C) = \Cliff(D'|_{C'})$, where the lattice $\Lambda = \langle D', \Delta' \rangle$ for effective divisors $D',\Delta'$ on $X'$.  Indeed, let $L_{D'}$ (respectively, $L_{\Delta'}$) be the deformation of the line bundle $\Oc(D)$ (respectively, $\Oc(\Delta)$) in the Noether--Lefschetz divisor $\NL_{\Lambda,g}$, and observe that any section of $\Oc(D)$ (respectively, $\Oc(\Delta)$) deforms to a section of $L_{D'}$ (respectively, $L_{\Delta'}$) as $H^1(X, \Oc(D)) = H^1(X,\Oc(\Delta)) = 0$.  Therefore, $L_{D'} = \Oc(D')$ and $L_{\Delta'} = \Oc(\Delta')$ for some effective divisors $D', \Delta'$ on $X'$.  These have the same intersection numbers as $D, \Delta$ by construction and hence $\Lambda = \langle D', \Delta'\rangle$.  Furthermore, because $H^1(X,D) = H^1(X,-\Delta)= 0$, by upper semicontinuity, we must have $H^1(X',D') = H^1(X',-\Delta') = 0$.  This forces $H^i(X', \Oc(D')) = H^i(X,\Oc(D))$ and $H^i(X', \Oc(-\Delta')) = H^i(X,\Oc(-\Delta))$.  From the exact sequence 
\[ 0 \to \Oc(-\Delta') \to \Oc(D') \to \Oc_{C'}(D'|_{C'}) \to 0\]
it follows that $\Cliff(D'|_{C'}) = \Cliff(D|_C)$. Assuming $C'$ is not exceptional, we have \[\gon(C') - 2 = \Cliff(C') \le \Cliff(D'|_{C'}) = \Cliff(D|_C) = \gon(C) - 2 \le \gon(C') - 2\] so each inequality must be an equality and hence $e= \gon(C) = \gon(C')$.  This implies $(X',L') \in \pgonlocus{e}{g}^{m}$, as desired. 

Now, assume $C'$ is exceptional for every smooth $C' \in |mL'|$.  From the chain of inequalities 
\[ \gon(C') - 3 = \Cliff(C') \le \Cliff(D'|_{C'}) = \Cliff(D|_C) = \gon(C) - 2 \]
and the fact that $\gon(C') \ge \gon(C) = e$, we conclude either $\gon(C') = e$ or $\gon(C') = e+1$.  If $\gon(C') = e$, then $(X',L') \in \pgonlocus{e}{g}^{m}$, as desired.  Now assume $\gon(C') = e+1$, and observe that this implies $\Cliff(C') = \Cliff(D'|_{C'}) = \Cliff(C)$.  Because $C'$ is exceptional, either $(X',L')$ is a Donagi--Morrison example or a generalized ELMS example. If it is a Donagi--Morrison example, we necessarily have $g = 2$, $m =3$, and $\gon(C') = 5$, so then $e = 4$. But then $|3L'|$ contains smooth curves of gonality 4 and in this case $\overline{\pgonlocus{4}{2}^{3}} = \Fc_2$.

Now, assume $C'$ is a generalized ELMS example.  We will use Lemma \ref{lem:ELMS-specialization} to show this is impossible by proving that the intersection matrix of $D$ and $\Delta$ is as in the Lemma, which would force $C$ itself to be exceptional, contradicting our assumption.  Because $\Cliff(C') = \Cliff(D'|_{C'})$, by the definition of a generalized ELMS example we see that $D'|_{C'}$ is the unique line bundle on $C'$ that computes $\Cliff(C')$, and hence $\Cliff(C') = \frac{g(C')-4}{2}$ by definition of the generalized ELMS examples.  However, $\Cliff(C') = e-2$ (where $\gon(C') = e+1$), which implies $e = \frac{g(C')}{2}$ or equivalently $e+1 = \frac{g(C')+2}{2}$ and hence all such $C'$ have the generic gonality.  Because $D'|_{C'}$ is the unique line bundle on $C'$ computing its Clifford index, for $D$ and $\Gamma$ as in the definition of the generalized ELMS examples, we claim we must have $D' = D + k \Gamma$ for some $k \in \Qb$.  Indeed, $D|_{C'}$ must agree with $D'|_{C'}$, and assuming that $X'$ is very general implies its Picard rank is exactly 2 and multiples of $\Gamma$ are the only classes that do not meet $C'$.  But, we have $(D')^2 = e -1$ by assumption and $D^2 = e-1$ by definition of the generalized ELMS examples, so 
\[ e-1 = (D+k \Gamma)^2 = D^2 + 2k - 2k^2 = e-1 + 2k(1-k) \] and therefore $k = 0$ or $k = 1$.  If $k = 0$, $D' \sim D$ so $C \sim D' + \Delta' \sim 2D + \Gamma$ so $\Delta' \sim D + \Gamma$.  If $k = 1$, then $D' \sim D + \Gamma$ and $C \sim D' + \Delta' \sim 2D + \Gamma$ implies $\Delta' \sim D$.  Up to switching $D'$ and $\Delta'$, we may assume $D' = D$ and $\Delta' = D + \Gamma$.  Then from the definition of the generalized ELMS examples, the intersection matrix of $D$ and $\Delta = D + \Gamma$ is of the form $\begin{pmatrix} e-1 & e \\ e & e-1 \end{pmatrix}$. By Lemma \ref{lem:ELMS-specialization}, the original surface $X$ and curve $C$ must have been a generalized ELMS example, so $C$ is exceptional, contradicting our assumption. Therefore, this case cannot occur. 

\textbf{Step 2(b):} To complete the proof, we consider the case that $C$ is exceptional. Then the membership of $(X, L) \in B_i$ for some $B_i$ will follow by separately analyzing each example classified by Knutsen (Theorem~\ref{thm:knutsen}).  In the case that $(X,L)$ is as in the Donagi--Morrison example, we have argued above that we must have $g = 2$, $m = 3$, and $e = 5$, and in this case the locus $\pgonlocus{5}{2}^{3}$ is dense in $\Fc_2$.
\par 
Now suppose that $(X, C)$ is a generalized ELMS example. Since $\Cliff(C)$ is computed uniquely by $\Oc_C(D)$ in this case, the rank $2$ sublattice $\Lambda = \left<D, D+\Gamma \right>$ is a Donagi--Morrison matrix and defines a Noether--Lefschetz divisor $ \NL_{\Lambda;g}$.  To complete the proof, it suffices to show the generic member of $\NL_{\Lambda;g}$ is contained in $\pgonlocus{e}{g}^{m}$ as this implies $\NL_{\Lambda;g} = B_i$ for some of the $B_i$ defined above.  For a very general point $(X', L') \in \NL_{\Lambda;g}$ the classes corresponding to the images of $D, \Gamma$ under $\Lambda \embed \Pic{X'}$ satisfy the same intersection relations as in a generalized ELMS example. Since the specialization homomorphism $\Pic{X'} \embed \Pic{X}$ is an injection and is compatible with the intersection product \cite[Chapter 17, Proposition 2.10]{Huybrechts16}, there does not exist a class $B \in \Pic{X'}$ so that $0 \le B^2 < D^2$ and $0 < B \cdot 
L'- B^2 \le D^2$. Therefore, $(X', L')$ is itself a generalized ELMS example and hence any smooth $C' \in |L'|$ is exceptional. Since $(X, C)$ and $(X', C')$ are both generalized ELMS examples for the same lattice, we see that $\gon(C) = \Cliff(C) + 3 = \Cliff(C') + 3 = \gon(C')$ and conclude that the generic point of $\NL_{\Lambda;g}$ is in the gonality locus $\pgonlocus{e}{g}^{m}$ as desired.
\end{proof}

\appendix 
\section{Examples of components of $\Gc_{e;g}$ in codimension $2$.}

In this section, we prove by example that the components of the unpolarized gonality loci $\gonlocus{e}{g}$ for $e \ge 3$ need not be divisorial.  The idea is similar in spirit to Theorem \ref{thm:codim2-hyperelliptic}.  

\begin{theorem}\label{thm:codim-2-comp}
Let $M$ be a $3 \times 3$ integer matrix and let $(g,h,e)$ be a triple of positive integers  satisfying the following properties:
\begin{enumerate}[label=(\alph*)]
    \item $M$ has signature $(1,2)$ and embeds primitively into the K3 lattice $E_8(-1)^{\oplus 2} \oplus U^{\oplus 3}$;
    \item the quadratic form defined by $M$ represents $2g - 2$ and $2 h - 2$ primitively;
    \item $e \ne 5$ is odd, $2< e < \floor{\frac{g+3}{2}}$, and $2<e < \floor{\frac{h+3}{2}}$; 
    \item the quadratic form defined by $M$ does not represent $0$, $ 2(2e - 3) - 2$, or $2(2e - 2) - 2$;
    \item $M$ contains a rank $2$ primitive Donagi--Morrison sub-lattice of type $(h,e)$ inside $\Fc_h$ which is a nonempty component of $\pgonlocus{e}{h} \subset \Fc_h$;
    \item For any integer $e'$ such that $2 < e' < e$, $M$ does not contain a rank $2$ Donagi--Morrison sub-lattice of type $(h,e')$ and, if $h = 2e$, $M$ does not contain the rank $2$ Donagi--Morrison sub-lattice of type $(h,e)$ given by $\begin{pmatrix} e-1 & e \\ e & e-1 \end{pmatrix}$;
    \item For any positive integer $h'$, 
    the only rank 2 Donagi--Morrison matrices of type $(h',e)$ that embed into $M$ admit only primitive embeddings into $M$ and do not represent $2g-2$.  
\end{enumerate}
Then, the locus $\gonlocus{e}{g} \subset \Fc_g$ has an irreducible component of codimension $2$. 
\end{theorem}

\begin{proof}
By property (a) and \cite[Chapter 14, Corollary 3.1]{Huybrechts16}, there exists a K3 surface $X$ with $\Pic{X} \cong \Zb^3$ and intersection form given by $M$. Let $L$ be the class representing $2g-2$, which is ensured by assumption (b). By reflecting along $(-2)$-classes as in Remark \ref{rem:reflections-nef}, we may assume that $L$ is nef. If $L$ was not basepoint-free, then by Proposition \ref{prop:base-locus-of-big-and-nef} there would exist a smooth elliptic curve $E \subset X$ such that $E^2 = 0$.  Since $M$ does not represent $0$ (assumption (d)), we conclude that the linear series $|L|$ is basepoint-free.  Now consider the class $C$ representing $2h-2$ from assumption (e).  By Lemma \ref{lem:DM-matrices-force-gon(C)}, $C$ may be assumed effective and either $|C|$ is not basepoint-free or $C$ may be assumed smooth.  In the first case, $M$ would represent $0$, a contradiction, so $|C|$ is basepoint-free.  From Lemma \ref{lem:DM-matrices-force-gon(C)}(2), we have $\gon(C) \le e$ if $C$ is not exceptional, or $\gon(C) \le e+1$ if $C$ is exceptional.  If $C$ is exceptional, by Theorem \ref{thm:knutsen}, either $X$ is a double cover of $\P^2$ branched over a plane sextic with $\Oc_X(C)$ the pullback of $\Oc_{\P^2}(3)$, which is not primitive, contradicting (a) and (b), or $X$ is a generalized ELMS example.  In this case, either $\gon(C) \le e$ or $ \gon(C) = e+1$ and by Example \ref{ex:ELMS}, $\Pic{X}$ contains a sub-lattice of the form $\begin{pmatrix} e-1 & e \\ e & e-1 \end{pmatrix}$ where by definition of the generalized ELMS examples, $h = 2e$, so this does not embed into $M$ by assumption (f).  In every case, this shows $\gon(C) \le e$.   
 
Furthermore, if $\gon(C) = e' < e$, then by assumption (c), either $C$ is hyperelliptic or Theorem \ref{thm:DonagiMorrison} applies.  We claim that both cases give a contradiction.  If $C$ is hyperelliptic, by \cite{Reid76a}, then the linear series $|C|$ realizes $X$ as a double cover of a Hirzebruch surface $\F_n$, $0 \le n \le 4$, or $\P^2$, the latter case being only possible if $h = 2$ or $h = 5$.  By Lemma \ref{lem:doublecover-lattice}, if $X$ is a cover of $\F_n$, the matrix $M$ represents $0$, contradicting assumption (d).  By assumption (c), $h > 2$, so the only remaining possibility is $h = 5$ and $X$ is a double cover of the Veronese surface.  However, in this case, $C$ is not primitive, contradicting assumption (b).   Therefore, $C$ is not hyperelliptic. However, because $C$ is not hyperelliptic, if $\gon(C) = e' < e$, there must be a sublattice for $(h,e')$ by Theorem \ref{thm:DonagiMorrison}. If this sublattice has rank 2 in $\Pic{X}$, then this contradicts assumption (f).  Otherwise, the classes $D, \Delta$ from Theorem \ref{thm:DonagiMorrison} are proportional, and as $C \in |D + \Delta|$ is smooth, $D + \Delta$ is big, nef, and basepoint-free.  Writing $D= n_1 P$ and $\Delta = n_2P$ for a primitive class $P$, it follows that $P$ is a quasi-polarization. By Lemma \ref{lem:prop-case}, we have $n_1 = 1$ and $C \in |(1+n_2)P|$, contradicting the primitivity of $C$. Therefore $\gon(C) = e$. 

Now suppose that $(X, L)$ was contained in a divisor $D \subset \Fc_g$ whose very general point lies in $\Gc_{e;g}$. First suppose $D$ is not a Noether--Lefschetz divisor. This would mean that a very general point $(X', L')$ of $D$ corresponds to a K3 surface of Picard rank $1$. By Proposition~\ref{prop:picard_rank_1_clifford}, the Clifford index of any smooth curve $C'$ on $X'$ is at least $\floor{(g-1)/2} > e - 2$ (using assumption (c)) so $\gon(C') > e$, a contradiction. 
\par 
To complete the proof, we need to show that the very general point of every Noether--Lefschetz divisor in $\Fc_g$ containing $(X, L)$ corresponds to a pair $(X', L')$ where $X'$ does not contain a smooth curve of gonality $e$.  Suppose for contradiction this is not the case and let $C' \subset X'$ be a smooth curve with gonality $e$.  As $e > 2$ by assumption, there are two possibilities: either $\gon(C')$ is the generic value of $\floor{\frac{g(C') + 3}{2}}$ or Theorem \ref{thm:DonagiMorrison} applies.  In the first case, $e = \gon(C') = \floor{\frac{g(C') + 3}{2}}$ so $g(C') \in \{ 2 e - 3, 2e - 2 \}$. This implies that $C'^2 = 2 g(C') - 2 \in \{ 2(2e - 3) - 2, 2 (2e - 2) - 2 \}$.  However, $\Pic{X'} \embed \Pic{X}$ by specialization giving a class representing either $2(2e - 3) - 2$ or $2(2e - 2) - 2$, violating assumption (d). 

In the other case, Theorem \ref{thm:DonagiMorrison} and Corollary \ref{cor:DMmatrices} apply.  Let $D', \Delta'$ be the classes from Theorem \ref{thm:DonagiMorrison}.  Because $C' \in |D' + \Delta'|$ is smooth, $D'+ \Delta'$ is big, nef, and basepoint-free.  We claim that $D', \Delta'$ cannot be proportional.  If they are proportional, write $D' = n_1 P'$ and $\Delta' = n_2 P'$ for a primitive class $P'$ on $X'$.  Because $D' + \Delta' = (n_1+n_2)P'$ is big and nef, $P'$ is a quasi-polarization on $X'$.  By Lemma \ref{lem:prop-case}, we have $n_1 = 1$ and, if $C'$ is not exceptional, then $e = \gon(C') = n_2 P'^2$.  As $P'^2$ is even, this contradicts assumption (c).  If $C'$ is exceptional, then $X'$ is a double cover of $\P^2$ branched over a sextic and $e = \gon(C') = 5$, contradicting assumption (c). We conclude the classes $D'$ and $\Delta'$ cannot be proportional.

Therefore, $\Pic{X'}$ must contain a rank $2$ Donagi--Morrison matrix for $(h',e) = (g(C'), e)$. The embedding $\Pic{X'} \embed \Pic{X}$ is compatible with the intersection pairing, so assuming the only embeddings of the given Donagi--Morrison matrix into the Picard lattice of $X$ are primitive, the embedding into $\Pic{X'}$ must also be primitive. Hence, the Picard lattice of $X'$ must be exactly the Donagi--Morrison matrix.  This matrix does not represent $2g-2$ by assumption (g), which contradicts the fact that $\Pic{X'}$ represents $2g-2$ since $X' \in \Fc_g$. Thus, we reach a contradiction in both cases and conclude that the locus $\Gc_{e;g} \subset \Fc_g$ is not purely divisorial and contains an irreducible component of codimension 2, explicitly given as the set of K3 surfaces of genus $g$ for which the Picard lattice contains $M$.
\end{proof}

\begin{remark}\label{rmk:assump(d)}
    Replacing assumption (d) in Theorem \ref{thm:codim-2-comp} with 
\begin{itemize}
\item[$(d')$] \textit{the quadratic form defined by $M$ does not represent $-2$, $2(2e-3) - 2$, $2(2e -2) - 2$, and $M$ does not contain the $2 \times 2$ sublattice $\begin{pmatrix} 0 & 2 \\ 2 & 0 \end{pmatrix}$},
\end{itemize}
then the same conclusion applies to find components of codimension 2 when $M$ is allowed to represent $0$. 
\end{remark}

\begin{example}
AI models were used to find the following matrices which satisfy the conditions of Theorem~\ref{thm:codim-2-comp}. 

One can verify that for $(g,h,e) = (15,13,7)$ the matrix
\[
    M=
    \begin{pmatrix}
        2 & 7 & 3\\
        7 & 8 & 5\\
        3 & 5 & -24
    \end{pmatrix}
\]
satisfies the conditions of the theorem above and hence it gives an example of a codimension 2 gonality 7 locus in $\Fc_{15}$.

Similarly, one can use Remark~\ref{rmk:assump(d)} for $(g,h,e) = (g,5,3)$. The matrix
\[
    M=
    \begin{pmatrix}
        0 & 3 & 3\\
        3 & 2 & 5\\
        3 & 5 & 2
    \end{pmatrix}
\]
gives an example of a gonality 3 locus in any $\Fc_{g}$ with $g > 0$ and $g \equiv 7 \pmod 9$.
\end{example}

\bibliographystyle{amsalpha}
\bibliography{biblio}

\vspace{0.5cm}

\footnotesize{
\noindent \textsc{University of Maryland} \\
\textit{E-mail address:} 
\href{mailto:dbejleri@umd.edu}{dbejleri@umd.edu} \\

\noindent \textsc{Instituto de Matem\'{a}tica Pura e Aplicada} \\
\textit{E-mail address:} \href{mailto:nathanchenmath@gmail.com}{nathanchenmath@gmail.com} \\

\noindent \textsc{Massachusetts Institute of Technology} \\
\textit{E-mail address:} \href{mailto:bvchurch@mit.edu}{bvchurch@mit.edu} \\

\noindent \textsc{University of California San Diego} \\
\textit{E-mail address:} \href{mailto:kedevleming@ucsd.edu}{kedevleming@ucsd.edu}
}

\end{document}